\documentclass[12pt,a4paper]{article}
\usepackage{amsthm}
\usepackage{amscd}
\usepackage{eucal}
\usepackage[active]{srcltx} 
\usepackage{latexsym}
\usepackage{amssymb}
\usepackage[centertags]{amsmath}
\usepackage{newlfont}
\usepackage{longtable}
\usepackage[doublespacing]{setspace}
\usepackage{graphicx}
\usepackage{graphics}
\usepackage{colortbl}
\usepackage{hyperref}
\usepackage{color}

\newtheorem{thm}{Theorem}
\newtheorem{lem}[thm]{Lemma}
\newtheorem{prop}[thm]{Proposition}

\numberwithin{equation}{section}

\newtheorem{Def}{Definition}
\newtheorem{rmk}{Remark}

\newcommand{\bege}[1]{\begin{equation}\label{#1}}
\newcommand{\eeq}{\end{equation}}

\def\QED{\hfill $\blacksquare$}

\newcommand{\x}{\mathbf{x}}
\newcommand{\y}{\mathbf{y}}
\newcommand{\uu}{\mathbf{u}}
\newcommand{\dd}{\mathbf{d}}
\newcommand{\vn}{\mathbf{v}}
\newcommand{\bn}{\mathbf{b}}

\newcommand{\tita}{\boldsymbol{\theta}}
\newcommand{\Btil}{\widetilde{\mathbf{B}}_n}
\newcommand{\Sigtil}{\widetilde{\mathbf{\Sigma}}_n}
\newcommand{\Bhat}{\widehat{\mathbf{B}}_n}

\newcommand{\Gahat}{\widehat{\mathbf{\Gamma}}_n}
\newcommand{\Ga}{\mathbf{\Gamma}}
\newcommand{\Gatil}{\widetilde{\mathbf{\Gamma}}_n}
\newcommand{\Sighat}{\widehat{\mathbf{\Sig}}_n}

\newcommand{\uuhat}{\widehat{\mathbf{u}}}
\newcommand{\z}{\mathbf{z}}
\newcommand{\B}{\mathbf{B}}
\newcommand{\Bo}{\mathbf{B}_0}
\newcommand{\A}{\mathbf{A}}
\newcommand{\Ahat}{\widehat{\mathbf{A}}}
\newcommand{\V}{\mathbf{V}}
\newcommand{\M}{\mathbf{M}}
\newcommand{\Vhat}{\mathbf{\widehat{V}}}
\newcommand{\Sig}{\mathbf{\Sigma}}
\newcommand{\U}{\mathbf{U}}

\newcommand{\R}{\mathbf{R}}
\newcommand{\Z}{\mathbf{Z}}
\newcommand{\T}{\mathbf{T}}
\newcommand{\Nor}[1]{N_q(\mathbf{0},#1)}
\newcommand{\Norpq}[1]{N_{qp}(\mathbf{0},#1)}
\newcommand{\Norma}[1]{N_{#1}(\mathbf{0},\mathbf{I})}
\newcommand{\Real}{\mathbb{R}}
\newcommand{\sighat}{\hat{\sigma}_n}
\newcommand{\Id}{\mathbf{I}}
\newcommand{\eps}{\varepsilon}
\newcommand{\ve}[1]{\text{vec}(#1)}

\begin{document}

\title{Estimates of MM type for the multivariate linear model}
\author{Nadia L. Kudraszow$^{\text{a}}$ and Ricardo A. Maronna$^{\text{b}}$\\
\small$^{\text{a}}$Universidad Nacional de La Plata\\
\small$^{\text{b}}$University of La Plata and C.I.C.P.B.A. (rmaronna@retina.ar)}
\date{}
\maketitle

\begin{abstract}
We propose a class of robust estimates for multivariate linear models. Based on the approach of MM estimation (Yohai 1987, \cite{Y}), we estimate the regression coefficients and the covariance matrix of the errors simultaneously. These estimates have both high breakdown point and high asymptotic efficiency under Gaussian errors. We prove consistency and asymptotic normality assuming errors with an elliptical distribution. We describe an iterative algorithm for the numerical calculation of these estimates. The advantages of the proposed estimates over their competitors are demonstrated through both simulated and real data. 

\end{abstract}


\section{Introduction}

\hspace{0.6cm}Consider a multivariate linear model (MLM) with random predictors, i.e., we observe $n$ independent identically distributed (i.i.d.) $(p+q)$-dimensional vectors, $\z_i=(\y_i',\x_i')$ with $1\leq i\leq n$, where
$\y_i=(y_{i1},\ldots ,y_{iq})'\in \Real^{q}, \;\;\; \x_i=(x_{i1},\ldots,x_{ip})'\in \Real^{p}$
and $'$ denotes the transpose. The $\y_i$ are the response vectors and the $\x_i$ are the predictors and both satisfy the equation
\begin{equation}\label{MLM}
\y_i=\Bo'\x_i+\uu_i \;\;\;\;\;\; 1\leq i\leq n,
\end{equation}
where $\Bo\in \Real^{p\times q}$ is the matrix of the regression parameters  and $\uu_i$ is a $q$-dimensional vector independent of $\x_i$.
If $\x_{ip}=1$ for all $ 1 \leq i \leq n $, we obtain a regression model with intercept.

We denote the distributions of $\x_i$ and $\uu_i$ by $G_0$ and $F_0$, respectively, and $\Sig_0$ is the covariance matrix of the $\uu_i$. The $p$-multivariate normal distribution with mean vector $\boldsymbol{\mu}$ and covariance matrix $\Sig$ is denoted by $N_p(\boldsymbol{\mu},\Sig)$.

In the case of $\uu_i$ with distribution $\Nor{\Sig_0}$,  the maximum likelihood estimate (MLE) of  $\Bo$ is the least squares estimate (LSE), and the MLE of $\Sig_0$ is the sample covariance matrix of the residuals. It is known that these estimates are not robust: a small fraction of outliers may have a large effect on their values.

Several approaches have been proposed to deal with this problem.
The first proposal of a robust estimate for the MLM was given Koenker and Portnoy \cite{Koenker Portnoy}. They proposed to apply a regression M-estimator, based on a convex loss function, to each coordinate of the response vector. The problems with this estimate is lack of affine equivariance and zero breakdown point. Several other estimates without these problem were defined later. Rousseeuw et al. \cite{Rousseeuw04} proposed estimates for the MLM based on a robust estimate of the covariance
matrix of $\z = (\x', \y')$. Bilodeau and Duchesne \cite{Bilodeau Duchesne} extended the S-estimates
introduced by Davies \cite{Davies} for multivariate
location and scatter; then Van Aelst and Willems \cite{Aelst Willems} studied the robustness of these estimators. Agull{\'o} et al. \cite{Agullo} extended the minimum covariance determinant
estimate introduced by Rousseeuw \cite{Rousseeuw85} and Roelandt et al. \cite{GS} extended the definition of GS-estimates introduced by Croux et al. \cite{Croux}. These estimates have a high breakdown point but are not highly efficient when the errors are Gaussian and $q$ is small. In order to solve this, Agull{\'o} et
al. \cite{Agullo} improved the efficiency of their estimates, maintaining their high breakdown point, by
considering one-step reweighting and one-step Newton-Raphson GM-estimates. Garc\'ia Ben et al. \cite{Tau multi} extended $\tau$-estimates for multivariate regression, obtaining a estimate with high breakdown point and a high Gaussian efficiency. Another important approach to obtain robust and efficient estimates is contrained M (CM)
estimation, proposed by Mendes and Tyler \cite{Mendes} for regression and by Kent and Tyler \cite{Kent}
for multivariate location and scatter. The bias of CM estimates for regression was studied by
Berrendero et al. \cite{Berrendero}. Following this approach, Bai et al. \cite{Bai} proposed CM estimates for
the multivariate linear model.

In this paper we propose robust estimates for the linear model based on the MM approach, first proposed by Yohai \cite{Y}  for the univariate linear
model, and later by Lopuha\"a \cite{Lopuaa}, Tatsuoka et al. \cite{Tatsuoka} and  Salibi\'an-Barrera et al. \cite{MM pos y esc} for multivariate location and
scatter. We show that our estimates have both a high breakdown point and a
high normal efficiency.

In Section 2 we define MM-estimates for the MLM and prove some properties. In Section 3 and 4 we study their breakdown point and Influence Function. In Section 5 and 6 we study the asymptotic properties
(consistency and asymptotic normality) of the MM-estimates assuming random predictors
and errors with an elliptical unimodal distribution. In Section 7 we describe a computing algorithm based on an iterative weighted MLE. In Section 8 we present the results of a simulation study and a real example in Section 9. All the proofs can be found in the Appendix.

\section{Definition and properties}

Before defining our class of robust estimates for the MLM, we will define a robust estimate of scale.

\begin{Def}\label{M-esc-def}Given a sample of size $n$, $\vn=(v_1,\ldots, v_n)$, an \textit{M-estimate of scale} $s(\vn)$ is defined as the value of $s$ that is solution of 
\begin{equation}\label{M-escala}
\displaystyle{\frac{1}{n}\displaystyle{\sum_{i=1}^n\rho_0\left(\frac{v_i}{s}\right)}=b,}
\end{equation}
where 
$b\in(0,1)$, or $s=0$ if $\sharp(v_i = 0) \geq n(1-b)$, where $\sharp$ is the symbol for cardinality.
\end{Def}

In this paper we use $b = 0.5$, which ensures the maximal asymptotic breakdown point (see \cite{Huber}).

The function $\rho_0$ should satisfy the following definition.

\begin{Def}\label{rhofuncion}
A \textit{$\rho$-function} will denote a function $\rho(u)$ which is a continuous nondecreasing function of $|u|$ such that $\rho(0)=0$, $\sup_u\rho(u)=1$, and $\rho(u)$ is increasing for nonnegative $u$ such that $\rho(u)<1$.
\end{Def}
Note that according to the terminology of Maronna et al. \cite{Maronna} this would be a ``bounded $\rho$-function''.
A popular $\rho$-function is the \textit{bisquare function}:
\begin{equation}\label{bicuadrada}
\rho_B(u)=1-(1-u^2)^3I(|u|\leq 1),
\end{equation}
where $I(\cdot)$ is the indicator function.

\begin{Def}\label{Mahalanobis}
Given a vector $\uu$ and a positive definite matrix $\mathbf{V}$, the \textit{Mahalanobis norm of $\uu$ with respect to $\mathbf{V}$} is defined as
\begin{equation*}
d(\uu,\V)=(\uu'\V^{-1}\uu)^{1/2}.
\end{equation*}

For particular given $\B\in\Real^{p\times q}$ and $\Sig\in \Real^{q\times q}$, we denote by $d_i(\B,\Sig)$ ($i=1,\dots, n$) the Mahalanobis norms of the residuals with respect to the
matrix $\Sig$, that is, 
\[d_i(\B,\Sig)=(\uuhat_i(\B)'\Sig^{-1}\uuhat_i(\B))^{1/2},\]
with  $\uuhat_i(\B)=\y_i-\B'\x_i$. 
\end{Def}

Using the concepts defined before, we can describe an \textit{MM-estimate for the MLM } by the following procedure:

\noindent Let $(\Btil,\Sigtil)$ be an initial estimate of $(\B_0,\Sig_0)$, with high breakdown point and such that $|\Sigtil|=1$, where  $|\Sigtil|$ is the determinant of $\Sigtil$ (i.e. $\Sigtil$ is an estimate of the shape of $\Sig_0$, $\Sig_0/|\Sig_0|^{1/q})$. Compute the Mahalanobis  norms of the residuals using $(\Btil,\Sigtil)$,
\begin{equation}\label{norma-resid}
d_i(\Btil,\Sigtil)=(\uuhat_i'(\Btil)\Sigtil^{-1}\uuhat_i(\Btil))^{1/2} \hspace{2em} 1\leq i\leq n.
\end{equation}
Then, compute the M-estimate of scale $\sighat:=s(\dd(\Btil,\Sigtil))$ of the above norms, defined by \eqref{M-escala}, using a  function $\rho_0$ as specified in Definition \ref{rhofuncion} and $b=0.5$. 

Let $\rho_1$ be another $\rho$-function such that
\begin{equation}\label{relacion-rhos}
\rho_1\leq \rho_0
\end{equation}
and let $\mathcal{S}_q$ be the set of all positive definite symmetric $q\times q$ matrices.

Let $(\Bhat,\Gahat)$ be any local minimum of
\begin{equation}\label{S}
	S(\B,\Ga)=\sum_{i=1}^n\rho_1\left(\frac{d_i(\B,\Ga)}{\sighat}\right)
\end{equation}
in $\Real^{p\times q}\times \mathcal{S}_q$, which satisfies
	\begin{equation}\label{MM3}
S(\Bhat,\Gahat)\leq S(\widetilde{\B}_n,\widetilde{\Sig}_n)
\end{equation}
and 	$|\Gahat|=1$. Then the MM-estimate of $\B_0$ is defined as $\Bhat$, and the respective estimate of $\Sig_0$ is
\begin{equation}\label{MM2}	
\Sighat=\sighat^2\Gahat.
\end{equation}

In the MM-estimates for the univariate linear model the residuals are used as a tool of outlier detection, in the MM-estimates for the multivariate linear model the Mahalanobis norms of the residuals play the same role. To compute the M-escale it is necessary to have an initial estimate of $\B_0$, to compute the residuals, and an initial estimate of the shape of $\Sig_0$ to compute the Mahalanobis norms of the residuals.

\begin{rmk}\label{cons-c}
One form of choosing the $\rho$-functions $\rho_0$ and $\rho_1$ in such a way that they satisfy \eqref{relacion-rhos} is the following. Let $\rho$ be a $\rho$-function and let $0<c_0<c_1$. We take 
\begin{equation}\label{c0c1}
\rho_0=\rho(u/c_0)\;\;\; \text{ and }\;\;\; \rho_1=\rho(u/c_1).
\end{equation}
The value $c_0$ should be chosen such that the asymptotic value of $\sighat$ is one when the errors $\uu_i$, with $i=1,\dots,n$, have distribution $\Norma{q}$. The choice of $c_1$ will determine the asymptotic efficiency of the MM-estimate. For more details see Remark \ref{ERA}.
\end{rmk}

The following theorem implies that the absolute minimum of $S(\B,\Ga/|\Ga|^{1/q})$ in $\Real^{p\times q}\times \mathcal{S}_q$ exists. Clearly, from this absolute minimum we can obtain an MM-estimate. However, any other local minimum $(\B,\Ga)$ which satisfies  \eqref{MM3}, may also be used to get an MM-estimate with high breakdown point and with high efficiency under Gaussian errors.   

Before stating the theorem we define $k_n$ as the maximum number of observations $(\y_i',\x_i')$ of a sample that are in a hyperplane, i.e., 
\begin{equation}\label{kn}	
	k_n:=\max_{\|\mathbf{v}\|+\|\mathbf{w}\|>0}\#\{i:\mathbf{v}'\x_i+\mathbf{w}'\y_i=\mathbf{0}\}.
\end{equation}

\begin{thm}\label{existencia}
Let $\Z=\{\z_1,\ldots, \z_n\}$ be a sample of size $n$ satisfying the MLM \eqref{MLM}, where $\z_i=(\y_i',\x_i')$.
If $k_n/n< 0.5$ then there is a pair $(\Bhat,\Gahat)$ that minimizes the function $S(\B,\Ga)$, defined in \eqref{S}, for all $(\B,\mathbf{\Gamma})\in \Real^{p\times q}\times \mathcal{S}_q$  such that  ${|\mathbf{\Gamma}|=1}$.
\end{thm}
The proof of this theorem can be found in the Appendix.

In the following theorem we obtain the estimating equations of MM-estimates.
\begin{thm}\label{teoMMec}
Assume that $\rho_1$ is differentiable. Then the MM-estimates $(\Bhat,\Sighat)$ satisfy the following equations:
\begin{equation}\label{MMec1}	
\displaystyle{\sum_{i=1}^nW\left(d_i(\Bhat,\Sighat)\right)\uuhat_i(\Bhat)\x_i'=\mathbf{0}}
\end{equation}
\begin{equation}\label{MMec2}	
\displaystyle{\Sighat=q\frac{\sum_{i=1}^n W\left(d_i(\Bhat,\Sighat)\right)\uuhat_i(\Bhat)\uuhat_i(\Bhat)'}{\sum_{i=1}^n \psi_1\left(d_i(\Bhat,\Sighat)\right)d_i(\Bhat,\Sighat)}}
\end{equation}
where $\psi_1(u)=\rho_1'(u)$ and $W(u)=\psi_1(u)/u$.
\end{thm}

\begin{rmk}\label{para algoritmo}
As we can see in equation \eqref{MMec1}, the $j$th column of $\Bhat$ is the weighted LSE  corresponding to the univariate regression whose dependent variable is the $j$th component of $\y$, the vector of independent variables  is the same that in the multivariate regression and the observation $i$ receives the weight $W\left(d_i(\Bhat,\Sighat)\right)$. Furthermore, by \eqref{MMec2}, $\Sighat$ is proportional to the sample covariance matrix of the weighted residuals with the same weights. As these weights depend on the estimates $\Bhat$ and $\Sighat$, we cannot use the relations \eqref{MMec1} and \eqref{MMec2} to compute the estimates, but they will be used to formulate an iterative algorithm in Section 6.
\end{rmk}

\begin{rmk}\label{equiv}
If  $\Btil$ is regression-, affine- and scale-equivariant and 
$\Sigtil$ is affine-equivariant and regression- and scale-invariant.
Then $\Bhat$ will be regression-, affine- and scale-equivariant and $\Sighat$ will be regression- and scale-invariant and affine-equivariant.
\end{rmk}
\section{Breakdown Point}

Now, to investigate the robustness of the MM-estimates, we will seek a lower bound of their finite sample breakdown point. The finite sample breakdown point of the coefficient matrix estimate is the smallest fraction of outliers that make the estimator unbounded, and the finite sample breakdown point of the covariance matrix estimate is the smallest fraction of outliers that make the estimate unbounded or singular.

Let $\Z=\{\z_1,\ldots, \z_n\}$ be a sample of size $n$ that satisfies the MLM (\ref{MLM}), where $\z_i=(\y_i',\x_i')$ and let $\widehat{\B}$ and $\widehat{\Sig}$ be estimates of $\B_0$ and $\Sig_0$ respectively.
We define
\[	
	\mathcal{Z}_m=\{\Z^*=\{\z_1^*,\ldots, \z_n^*\}\text{ such that } \sharp\{i:\z_i=\z_i^*\}\geq n-m\},
\]
\[
S_m(\Z,\widehat{\B})=\sup\{\|\widehat{\B}(\Z^*)\|_2 \text{ with }\Z^*\in\mathcal{Z}_m\},
\]
\[
S_m^+(\Z,\widehat{\Sig})=\sup\{\lambda_1(\widehat{\Sig}(\Z^*)) \text{ with }\Z^*\in\mathcal{Z}_m\}
\]
y 
\[
S_m^-(\Z,\widehat{\Sig})=\inf\{\lambda_q(\widehat{\Sig}(\Z^*)) \text{ with }\Z^*\in\mathcal{Z}_m\},
\]
where $\lambda_1(\widehat{\Sig}(\Z^*))$ and $\lambda_q(\widehat{\Sig}(\Z^*))$ are the  largest and smallest eigenvalue of $\widehat{\Sig}(\Z^*)$ respectively.

\begin{Def}
The \textit{finite sample breakdown point of $\widehat{\B}$} is $\varepsilon^*(\Z,\widehat{\B})=m^*/n$ where 
\[m^*=\min\{m:S_m(\Z,\widehat{\B})=\infty\},
\]
the \textit{finite sample breakdown point of $\widehat{\Sig}$} is $\varepsilon^*(\Z,\widehat{\Sig})=m^*/n$ where  
\[
m^*=\min\{m:\frac{1}{S_m^-(\Z,\widehat{\Sig})} + S_m^+(\Z,\widehat{\Sig})=\infty\}
\]
and $\varepsilon^*_n(\Z,\widehat{\B},\widehat{\Sig})=\min\{\varepsilon^*(\Z,\widehat{\B}),\varepsilon^*(\Z,\widehat{\Sig})\}$.
\end{Def}

The following theorem gives a lower bound for the breakdown point of MM-estimates.

\begin{thm}\label{TeoBP} Let $\Z=\{\z_1,\ldots, \z_n\}$, with $\z_i=(\y_i',\x_i')$ that satisfies the MLM (\ref{MLM}), and $k_n$ defined in \eqref{kn}.  Consider $\rho_0$ and $\rho_1$ two $\rho$-functions that satisfy \eqref{relacion-rhos} and suppose that $k_n<n/2$.
Then  
\begin{equation}\label{BP}	
\varepsilon^*_n(\Z,\Bhat,\Sighat)\geq\min\left(\varepsilon^*_n(\Z,\Btil,\Sigtil),\frac{[n/2]-k_n}{n}\right).
\end{equation}
\end{thm}

Since $k_n$ is always greater or equal than $p+q-1$, if $\varepsilon^*_n(\Z,\Btil,\Sigtil)$ is close to $0.5$ the maximum lower bound will be $({[n/2]-(p+q-1)})/{n}$, i.e. when the points are in general position the finite sample breakdown point is close to $0.5$ for large $n$.

If we didn't fix $b=0.5$ and  if $k_n<n(1-b)$, we would have the same bound as in \eqref{BP} but with $[n(1-b)]$ in place of $[n/2]$.
In this case, the maximum finite sample breakdown point would be attained in $b=0.5-k_n/n$ which is very close to our choice of $b=0.5$ when $k_n/n$ is small. 

%


\section{Influence function}

Consider an estimate $\widehat{\boldsymbol{\tita}}_n$ depending on a sample $\Z = \{\z_1,\dots , \z_n\}$ of i.i.d. variables in $\Real^k$ with distribution $H_{\boldsymbol{\tita}}$, where $\boldsymbol{\tita}\in\Theta\subset\Real^m$. Let $T$ be an estimating functional  of $\tita$ such that $\T(H_n)=\widehat{\boldsymbol{\tita}}_n$, where $H_n$ is the corresponding empirical distribution. Suppose that $\T$ is Fisher consistent, i.e. $\T(H_{\boldsymbol{\tita}})=\tita$. The \textit{influence function} of $\T$, introduced by Hampel \cite{Hampel}, measures the effect on the functional of a small fraction of point mass contamination. If $\delta_\z$ denotes the probability distribution that assigns mass 1 to $\x$, then the influence function  is defined by
\[
IF(\z,\T,{\boldsymbol{\tita}})=\displaystyle{\lim_{\varepsilon \to 0}}\frac{\T((1-\varepsilon)H_{\boldsymbol{\tita}}+\varepsilon\delta_{\z})-\T(H_{\boldsymbol{\tita}})}{\varepsilon}=\displaystyle{\left.\frac{\partial \T((1-\varepsilon)H_{\boldsymbol{\tita}}+\varepsilon\delta_{\z})}{\partial\varepsilon}\right|_{\varepsilon=0}},
\]

In our case, $\z=(\y',\x')'$ satisfy the linear model $\eqref{MLM}$, $\boldsymbol{\tita}=(\B_0,\Sig_0)$ and $H_{\boldsymbol{\tita}}=H_0$. Let $\T_{0,1}$, $\T_{0,2}$ be the functional estimates asociated to the inicial estimates $\Btil$ and $\Sigtil$, and $\T_1$, $\T_2$ the functional estimates  corresponding to the MM-estimates $\Bhat$ and $\Sighat$. Then, according to \eqref{MMec1} and \eqref{MMec2}, given a distribution function $H$ of $(\y',\x')'$, the pair $(\T_1(H),\T_2(H))$ is the value of $(\B,\Sig)$ satisfying
\begin{equation*}\label{FT1}	
\displaystyle{E_{H}W\left(d(\B,\Sig)\right)\uuhat(\B)\x'=\mathbf{0}},
\end{equation*}
\begin{equation*}\label{FT2}	
\displaystyle{\Sig=q\frac{E_{H} W\left(d(\B,\Sig)\right)\uuhat(\B)\uuhat(\B)'}{E_{H_0} \psi_1\left(d(\B,\Sig)\right)d(\B,\Sig)}},
\end{equation*}
and
\[
\Sig=S(H)^2\Ga, \text{ with } |\Ga|=1,
\]
where $d(\B,\Sig)=d(\uuhat(\B),\Sig)$, $\uuhat(\B)=\y-\B'\x$ and 
\[
\displaystyle{E_{H}\rho_0\left(\frac{d(\T_{0,1}(H),\T_{0,2}(H))}{S(H)}\right)=\mathbf{0}}.
\]

Note that the M-estimate of scale, $\sighat$, used in the definition of MM-estimates $(\Bhat,\Sighat)$, verify $\sighat=S(H_n)$, where $H_n$ is the empirical distribution  of $\z_1,\dots,\z_n$.

Next we will state the influence function of MM-estimators for the case where errors in \eqref{MLM} have an elliptical distribution with unimodal density. For that, we need to make the following assumptions:\\

\noindent(A1) $\rho_1$ is strictly increasing in $[0,\kappa]$ and constant in $[\kappa,+\infty)$ for some constant $\kappa<\infty$.

\noindent (A2) $P_{G_0}(\B'\x=0)<0.5$ for all $\B\in\Real^{p\times q}$.

\noindent(A3) The distribution $F_0$ of $\uu_i$ has a density of the form 

\begin{equation}\label{densidad}
f_0(\uu)=\frac{f_0^*(\uu'\Sig_0^{-1}\uu)}{|\Sig_0|^{1/2}}
\end{equation}

\noindent where $f_0^*$ is nonincreasing and has at least one point of decrease in the interval where $\rho_1$ is strictly increasing.

\noindent(A4) $G_0$ has second moments and $E_{G_0}(\x{\x}')$ is no singular.

\begin{thm}\label{TFI} Let $(\y',\x')$ be a random vector satisfying the MLM  (\ref{MLM}) with parameters $\B_0$ and $\Sig_0$. Assume that (A1)-(A4) hold and that the partial derivatives of $\displaystyle{E_{H}W\left(d(\B,\Sig)/S(H)\right)\uuhat(\B)\x'}$ can be obtained differentiating with respect to each parameter inside the expectation, where $H_0$ is the distribution of $(\y',\x')'$. Suppose that the functional estimates associated to the initial estimates $\Btil$ and $\Sigtil$ are affine-equivariant. Then, the influence function for the functional estimator $\T_1$ corresponding to the MM-estimate $\Bhat$ is
\begin{eqnarray*}
&&\hspace{-2em}IF(\z_0,\T_1,\B_0,\Sig_0)\\
&=& c^{-1}W\left(\frac{\left((\y_0-\B_0'\x_0)'\Ga_0^{-1}(\y_0-\B_0'\x_0)\right)^{1/2}}{\sigma}\right)E_{G_0}(\x\x')^{-1}\x_0(\y_0-\B_0'\x_0)',
\end{eqnarray*}
where $\z_0=(\y_0',\x_0')'\in \Real^{q+p}$, $\sigma=S(H_0)$, $\Ga_0=\Sig_0|\Sig_0|^{-1/q}$,\hspace{0.3cm} and 
\[
c=\frac{E_{F_0}W'\left((\uu'\Ga_0^{-1}\uu)^{1/2}/{\sigma}\right)(\uu'\Ga_0^{-1}\uu)^{1/2}}{q\sigma}+E_{F_0}W\left(\frac{(\uu'\Ga_0^{-1}\uu)^{1/2}}{\sigma}\right).
\]
\end{thm}

As in the case of MM-estimators for univariate linear regression, the influence function of the proposed MM-estimate is unbounded.

\section{Consistency}

We will now show the consistency of MM-estimates for multivariate regression for the case in which errors in (\ref{MLM}) have an elliptical distribution with an unimodal density. For this, we need the following additional assumptions:

\begin{thm}\label{cons} Let $(\y_i',\x_i')$, $1\leq i\leq n$, be a random sample of the MLM  (\ref{MLM}) with parameters $\B_0$ and $\Sig_0$. Assume that $\rho_0$ and $\rho_1$ are $\rho$-functions that satisfy the relation \eqref{relacion-rhos}, that (A1)-(A3) hold and that the initial estimates $\Btil$ and $\Sigtil$ are consistent for $\B_0$ and $\Ga_0$ respectively, where $\Ga_0=\Sig_0|\Sig_0|^{-1/q}$; then the MM-estimates $\Bhat$ and $\Sighat$ satisfy

\begin{description}
	\item (a) $\lim_{n\rightarrow\infty} \Bhat = \B_0$ a.s..
	\item (b) $\lim_{n\rightarrow\infty} \Sighat = \sigma_0^2\Sig_0$ a.s.
with	$\sigma_0$ defined by 
\end{description}
\begin{equation}\label{4.2}
E_{F_0}\left(\rho_0 \left(\frac{\|\uu\|}{\sigma_0}\right)\right)=b.
\end{equation}
\end{thm}

\section{Asymptotic Normality}

Before obtaining the limit distribution of $\Bhat$ we need to make some additional assumptions.

\noindent(A5) $\rho_1$ is differentiable, $\psi_1=\rho_1'$ and $W(u)=\psi_1(u)/u$ is differentiable with bounded derivative. 

\noindent(A6) $E_{G_0}\|\x\|^4<\infty$, $E_{G_0}\|\x\|^6<\infty$, $E_{H_0}\|\x\|^4\|\y\|^2<\infty$ and $E_{H_0}\|\x\|^2\|\y\|^4<\infty$, where $H_0$ is the distribution of $\z=(\y',\x')'$.

\noindent (A7) Let $\tita=(\B,\Sig)$ and 
\begin{equation}\label{phi}
\phi(\z;\tita)=W\left(d(\B,\Sig)\right)\ve{(\y-\B'\x)\x'}.
\end{equation}
The function $\Phi(\tita)=E_{H_0}\phi(\z;\tita)$ has a partial derivative $\partial\Phi/\partial \ve{\B'}'$ which is continuous at $\tita_0=(\B_0,\sigma_0^2\Sig_0)$ and the matrix
\begin{equation}\label{Lambda}
\mathbf{\Lambda}=\frac{\partial\Phi(\B,\Sig)}{\partial \ve{\B'}'}{(\B_0,\sigma_0^2\Sig_0)}
\end{equation}
is non singular.

\begin{thm}\label{normalidad}
Let $\z_i=(\y_i',\x_i')$, with $1\leq i\leq n$, be a random sample from the model (\ref{MLM}) with parameters $\B_0$ and $\Sig_0$. Assume that the $\rho$-function  $\rho_1$ satisfies (A1), that (A3)-(A7) hold and  that the estimates $\Btil$ and $\Sigtil$ are consistent for $\B_0$ and $\Ga_0=\Sig_0|\Sig_0|^{-1/q}$ respectively; then $n^{1/2}\ve{\Bhat'-\B_0'}\stackrel{d}{\rightarrow}\Norpq{\V}$, where $\stackrel{d}{\rightarrow}$ denotes convergence in distribution and
\begin{equation}
\V=\mathbf{\Lambda}^{-1}\M{\mathbf{\Lambda}^{-1}}'
\end{equation}
where $\M$ is the covariance matrix $\phi(\z_1,(\B_0,\sigma_0^2\Sig_0))$, with $\phi$ defined in \eqref{phi} and $\mathbf{\Lambda}$ is defined in \eqref{Lambda}. 
\end{thm}

Assumptions (A4)-(A7) are sufficient to prove Theorem \ref{normalidad}, but we conjecture that the limit distribution of $\Bhat$ can be proved under less restrictive hypotheses.

\begin{rmk} Note that the rate of convergence of the MM-estimates depends only on the consistency but not on the rate of convergence, of the initial estimates.
\end{rmk}

Under suitable differentiability conditions we can obtain a more detailed expression of the covariance matrix $\V$ of Theorem \ref{normalidad}. 

\begin{prop}\label{V} 
If $W_1(u)=W(\sqrt{u})$ is differentiable with bounded derivative and the initial estimates $(\Btil,\Sigtil)$ are affine-equivariant, then
\begin{equation}\label{V explicita}
\displaystyle{\V=\left[\frac{\sigma_0^2}{q}E_{F_0}\left(\psi_1\left(\frac{v}{\sigma_0}\right)\right)^2\biggl/\left(E_{F_0}W^*\left(\frac{v}{\sigma_0}\right)\right)^{2}\right](E_{G_0}\x\x')^{-1}\otimes\Sig_0}
\end{equation}
where
\begin{equation}\label{Westrella*}
W^*\left(\frac{v}{\sigma_0}\right)=\frac{2}{q\sigma_0^2}W_1'\left(\frac{v^2}{\sigma_0^2}\right)v^2+W\left(\frac{v}{\sigma_0}\right)
\end{equation}
and
\[
v=\left(\uu'\Sig_0^{-1}\uu\right)^{1/2}.
\]
\end{prop}

From the proof of Proposition \ref{V} (see Appendix), it is easily seen that if $W_1(u)$ is continuously differentiable with bounded derivative, assumption (A7) holds if and only if $E_{F_0}W^*\left({\left(\uu'\Sig_0^{-1}\uu\right)^{1/2}}/{\sigma_0}\right)\neq0.$

\begin{rmk}\label{ERA}
The covariance matrix of the MLE, is given by 
\[\V=\left({E_{F_0}(v^2)}/{q}\right)(E_{G_0}\x\x')^{-1}\otimes\Sig_0.\] Then the asymptotic relative efficiency of the MM-estimate $\Bhat$ with respect to the MLE is
\begin{equation}\label{ERAformula}
ARE(\psi_1,F_0)=E_{F_0}(v^2)\frac{\left(E_{F_0}W^*\left(\frac{v}{\sigma_0}\right)\right)^{2}}{{\sigma_0^2}E_{F_0}\left(\psi_1\left(\frac{v}{\sigma_0}\right)\right)^2}.
\end{equation}
 
As we mentioned in Remark \ref{cons-c}, to obtain an MM-estimate which simultaneously has  high breakdown point and high efficiency under normal errors  it suffices to choose $c_0$ and $c_1$ in \eqref{c0c1} appropriately. The constant $c_0$ can be chosen so that
\begin{equation}\label{rel-c0}
E\left(\rho \left(\frac{(\uu'\Ga_0^{-1}\uu)^{1/2}}{c_0}\right)\right)=b,
\end{equation}
where $\uu$ is $\Nor{\Sig_0}$, $\Sig_0=|\Sig_0|^{1/q}\Ga_0$ y $b=0.5$, this ensures a high breakdown point and that the asymptotic relative efficiency \eqref{ERAformula} depends only of $c_1$. Then, $c_1$ can be chosen so that the MM-estimate has the desired efficiency without affecting the breakpoint that depends only of $c_0$.
\end{rmk}

Table \ref{tabla1} gives the values of $c_0$ verifying \eqref{rel-c0} for different values of $q$. Table \ref{tabla2} gives the values of $c_1$ needed to attain different levels of asymptotic efficiency. In both cases the function $\rho$ from \eqref{c0c1} is equal to the bisquare function,  $\rho_B$,  given in \eqref{bicuadrada}.

\begin{table}[h!]
\begin{center}\begin{tabular}{l|llllll} 
\hline
 \rowcolor[rgb]{0.8,0.8,0.8}  $q$\hspace{1.2cm}  & 1\hspace{1.2cm} &  2\hspace{1.2cm} &  3\hspace{1.2cm} & 4\hspace{1.2cm} &  5\hspace{1.2cm} & 10\\ 
\hline 
$c_0$ \hspace{0.5cm}& 1.56 &  2.66 & 3.45 & 4.10 & 4.65 & 6.77 \\ 
\hline 
\end{tabular}
\end{center}
\caption{Values of $c_0$ for the bisquare function.}
\label{tabla1}
\end{table}

\begin{table}[h!]
\begin{center}\begin{tabular}{l|llllll} 
\hline
\rowcolor[rgb]{0.8,0.8,0.8} $ARE$\hspace{0.5cm} &  \multicolumn{6}{>{\columncolor[rgb]{0.8,0.8,0.8}}c}{$q$  } \\ 
 \hline
 \rowcolor[rgb]{0.8,0.8,0.8} \hspace{1.2cm} & 1\hspace{1.2cm} &  2\hspace{1.2cm} &  3\hspace{1.2cm} & 4\hspace{1.2cm} &  5\hspace{1.2cm} & 10\\ 
\hline 
 0.80\hspace{0.5cm} &  3.14 &  3.51 &  3.82 & 4.10 & 4.34 & 5.39\\ 
 0.90\hspace{0.5cm} & 3.88 & 	4.28 &	4.62 &	4.91 &	5.18 &	6.38\\	
 0.95\hspace{0.5cm} & 4.68 &	5.12 &	5.48 &	5.76 &	6.10 &	7.67\\
\hline 
\end{tabular}
\end{center}
\caption{Values of $c_1$ for the bisquare function to attain given values of the asymptotic relative efficiency ($ARE$) under normal errors.}
\label{tabla2}
\end{table}


\section{Computing algorithm}

In this section we propose an iterative algorithm to compute $\Bhat$ and $\Sighat$ based on the Remark \ref{para algoritmo}. Let $\z_i=(\y_i',\x_i')$ be a sample of size $n$ and assume we have computed the initial estimates $\widetilde{\mathbf{B}}_{n}$ and $\widetilde{\mathbf{\Sig}}_{n}$ with high breakdown point and such that $|\widetilde{\mathbf{\Sig}}_{n}|=1$.

\begin{enumerate}
\item Using the initial values $\widetilde{\mathbf{B}}^{(0)}=\widetilde{\mathbf{B}}_{n}$ and $\widetilde{\mathbf{\Gamma}}^{(0)}=\widetilde{\mathbf{\Sig}}_{n}$, compute the M-estimate of scale $\sighat:=s(\dd(\widetilde{\mathbf{B}}^{(0)},\widetilde{\mathbf{\Gamma}}^{(0)}))$, defined by \eqref{M-escala}, using a function $\rho_0$ as in the definition and $b=0.5$ and the matrix $\widetilde{\mathbf{\Sigma}}^{(0)}=\hat{\sigma}_{n}^2\widetilde{\mathbf{\Gamma}}^{(0)}$.
\item Compute the weights $\omega_{i0}=W\left(d_i(\widetilde{\mathbf{B}}^{(0)},\widetilde{\mathbf{\Sigma}}^{(0)})\right)$ for $1\leq i\leq n.$ These weights are used to compute each column of $\widetilde{\mathbf{B}}^{(1)}$ separately by weighted least squares. 
\item Compute the matrix
\[
\displaystyle{{\widetilde{\mathbf{C}}^{(1)}}=\sum_{i=1}^n\omega_{i0}\uuhat_i(\widetilde{\mathbf{B}}^{(1)})\uuhat_i'(\widetilde{\mathbf{B}}^{(1)})},
\]
and with it the matrix $\widetilde{\mathbf{\Sigma}}^{(1)}=\hat{\sigma}_{n}^2\widetilde{\mathbf{C}}^{(1)}/|\widetilde{\mathbf{C}}^{(1)}|^{1/q}$.
\item Suppose that we have already computed $\widetilde{\mathbf{B}}^{(k-1)}$ and $\widetilde{\mathbf{\Sigma}}^{(k-1)}$. Then {$\widetilde{\mathbf{B}}^{(k)}$} and {$\widetilde{\mathbf{\Sigma}}^{(k)}$} are computed using the steps 2 and 3 but starting from $\widetilde{\mathbf{B}}^{(k-1)}$ and $\widetilde{\mathbf{\Sigma}}^{(k-1)}$ instead of $\widetilde{\mathbf{B}}^{(0)}$ and $\widetilde{\mathbf{\Sigma}}^{(0)}$.
\item The procedure is stopped at step $k$ if the relative absolute differences of all elements of the matrices $\widetilde{\mathbf{B}}^{(k)}$ and $\widetilde{\mathbf{B}}^{(k-1)}$ and the relative absolute differences of all the Mahalanobis norms of residuals $\uuhat_i(\widetilde{\mathbf{B}}^{(k)})$ and $\uuhat_i(\widetilde{\mathbf{B}}^{(k-1)})$ with respect to $\widetilde{\Sig}^{(k)}$ and $\widetilde{\Sig}^{(k-1)}$ respectively are smaller than a given value $\delta$.
\end{enumerate}

The following theorem, whose proof can be found in the Appendix, shows that the iterative procedure to compute MM-estimates yields the descent of the objective function. 

\begin{thm}\label{conv-algoritmo}
If $W(u)$ is nonincreasing in $|u|$ then at each iteration of the algorithm the function $\sum_{i=1}^n\rho_1\left(d_i(\mathbf{{B}},\mathbf{{\Sig}})\right)$ is nonincreasing.
\end{thm}

\section{Simulation}
\subsection{Simulation design}
To investigate the performance of the proposed estimates we performed a simulation study.\\
\noindent- We consider the MLM given by \eqref{MLM} for two cases: $p = 2$, $q = 2$ and $p = 2$,
$q = 5$. Due to the equivariance of the estimators we take, without loss of generality, $\B_0 = \mathbf{0}$ and
$\Sig_0 = \Id_q$. The errors $\uu_i$ are generated from an $\Norma{q}$ distribution and the predictors $\x_i$ from an $\Norma{p}$ distribution. 

\noindent- The sample size is $100$ and the number of replications is $1000$. We consider uncontaminated samples and samples that contain $10\%$ of identical outliers of the form $(\x_0, \y_0)$ with $\x_0 = (x_0, 0, \dots , 0)$ and $\y_0 = (mx_0, 0,\dots, 0)$.
The values of $x_0$ considered are $1$ (low leverage outliers) and $10$ (high leverage outliers). We take
a grid of values of $m$, starting at $0$. The grid was chosen in order that all robust estimates attain the maximum values of their error measure.\\
\noindent- Let $\widehat{\B}^{(k)}$ be the estimate of $\B_0$ obtained in the $k$th replication. Then, since we are taking $\B_0 = \mathbf{0}$, the estimate of the mean squared error (MSE) is given by
\[
\text{MSE}=\frac{1}{1000}\left(\sum_{k=1}^{1000}\sum_{i=1}^{p}\sum_{j=1}^{q}\left(\widetilde{\B}^{(k)}_{ij}\right)^2\right).
\]

It must be recalled that the distributions of robust estimates under contamination are themselves heavy-tailed, and it is therefore prudent to evaluate their performance through robust measures (see \cite{Huber} Sec. 1.4, p. 12, and \cite{Huber64} p.75). For this reason, we employed both MSE, and trimmed mean squared error (TMSE), which compute the $10\%$ (upper)
trimmed average of \[\left\{\sum_{i=1}^{p}\sum_{j=1}^{q}\left(\widetilde{\B}^{(k)}_{ij}\right)^2\right\}_{k=1}^{1000}.\]
 The results given below correspond to this MSE,
although the TMSE yields qualitatively similar results (in the uncontaminated case the results are the same).\\

\subsection{Description of the estimators}

For each case, four estimates are computed: the MLE, an S-estimate, a $\tau$-estimate and an MM-estimate. 

For the MLM, the S-estimates are defined by 
\[
(\widehat{\B},\widehat{\Sig})=\arg\min\{|\Sig|: (\B,\Sig)\in\Real^{p\times q}\times\mathcal{S}_q\}
\]
subject to 
\[
s^2(d_1(\B,\Sig),\dots,d_n(\B,\Sig))=q,
\]
where $s$ is an M-estimate of scale.

Garc\'ia Ben et al. \cite{Tau multi} extended $\tau$-estimates to the MLM by defining
\[
(\widehat{\B},\widehat{\Sig})=\arg\min\{|\Sig|: (\B,\Sig)\in\Real^{p\times q}\times\mathcal{S}_q\}
\]
subject to
\begin{equation}\label{kapa}
\tau^2(d_1(\B,\Sig),\dots,d_n(\B,\Sig))=\kappa,
\end{equation}
where the $\tau$-scale is defined by 
\begin{equation}\label{tauscale}
\tau^2(\vn)=(s^2(\vn)/{n})\sum_{i=1}^{n}\rho_2\left({|v_i|}/{s(\vn)}\right),
\end{equation}
where $\vn=(v_1,\dots,v_n)$, $\rho$ is a $\rho$-function and $s$ is an M-estimate of scale.

The robust estimates are based on bisquare $\rho$-functions. The M-estimate of scale used in the S-estimate is defined by $\rho_0(u) = \rho_B(u/c_0)$, and $b=0.5$ so that the S-estimate has
breakdown point $0.5$ (see Table \ref{tabla1}). The $\tau$-estimate uses the same $\rho_0$ and $b$ as the S-estimate to compute the M-scale and $\rho_2(u) = \rho_{B}(u/c_2)$, where $c_2$ is chosen together with the constant  $\kappa$, from equation \eqref{kapa}, so that the $\tau$-estimate has an ARE equal to $0.90$ when the errors are Gaussian (see Table 2 in \cite{Tau multi} in which $\kappa=6\kappa_2/c_2^2$). The initial estimate needed to compute the $\tau$-estimate is computed using 2000 subsamples. 
The MM-estimate uses the same  $\rho_0$ as the S-estimate to compute the M-estimate of scale and $\rho_1(u) = \rho_{B}(u/c_1)$, where $c_1$ is chosen so that the MM-estimate has an ARE equal to $0.90$ when the errors are Gaussian (see Table \ref{tabla2}). We use the S-estimates as $(\Btil,\Sigtil)$ and the value of $\delta$ in step 5 of the computing algorithm is taken equal to $10^{-4}$.

\begin{table}[t]
\small
\begin{center}\begin{tabular}{>{\columncolor[rgb]{0.8,0.8,0.8}}l|lllllllll}
\hline
 Estimate &  \multicolumn{4}{>{\columncolor[rgb]{0.8,0.8,0.8}}l}{\hspace{-0.2em}$q=2$}& \multicolumn{1}{>{\columncolor[rgb]{0.8,0.8,0.8}}l}{ \hspace{0em}}& \multicolumn{4}{>{\columncolor[rgb]{0.8,0.8,0.8}}l}{$q=5$}\\ 
  &  \multicolumn{1}{>{\hspace{-0.2em}\columncolor[rgb]{0.8,0.8,0.8}}l}{MSE}&  \multicolumn{1}{>{\columncolor[rgb]{0.8,0.8,0.8}}l}{\hspace{-0.2em}SE}&  \multicolumn{1}{>{\columncolor[rgb]{0.8,0.8,0.8}}l}{\hspace{-0.2em}REFF}&  \multicolumn{1}{>{\columncolor[rgb]{0.8,0.8,0.8}}l}{\hspace{-0.2em}ARE}&\multicolumn{1}{>{\columncolor[rgb]{0.8,0.8,0.8}}l}{\hspace{0em}}&\multicolumn{1}{>{\columncolor[rgb]{0.8,0.8,0.8}}l}{\hspace{-0.2em}MSE}&  \multicolumn{1}{>{\columncolor[rgb]{0.8,0.8,0.8}}l}{\hspace{-0.2em}SE}&  \multicolumn{1}{>{\columncolor[rgb]{0.8,0.8,0.8}}l}{\hspace{-0.2em}REFF}&  \multicolumn{1}{>{\columncolor[rgb]{0.8,0.8,0.8}}l}{\hspace{-0.2em}ARE} \\ 
\hline 
 MLE	&	\hspace{-0.2em}0.041 &	\hspace{-0.2em}0.001	& \hspace{-0.2em}1.00 &	\hspace{-0.2em}1.00 &\hspace{0em}& \hspace{-0.2em}0.103 &	\hspace{-0.2em}0.002 &	\hspace{-0.2em}1.00 &	\hspace{-0.2em}1.00\\
S-estimate	&	\hspace{-0.2em}0.074 &	\hspace{-0.2em}0.002 &	\hspace{-0.2em}0.55 &	\hspace{-0.2em}0.58 &\hspace{0em}& \hspace{-0.2em}0.125 &	\hspace{-0.2em}0.002 &	\hspace{-0.2em}0.83 &	\hspace{-0.2em}0.85\\
$\tau$-estimate & 		\hspace{-0.2em}0.046 &	\hspace{-0.2em}0.001 &	\hspace{-0.2em}0.89 &	\hspace{-0.2em}0.90 &\hspace{0em}& \hspace{-0.2em}0.116 &	\hspace{-0.2em}0.002 &	\hspace{-0.2em}0.90 &	\hspace{-0.2em}0.90\\
MM-estimate &		\hspace{-0.2em}0.046	 &\hspace{-0.2em}0.001 &	\hspace{-0.2em}0.89 &	\hspace{-0.2em}0.90 &\hspace{0em}& \hspace{-0.2em}0.116 &	\hspace{-0.2em}0.002 &	\hspace{-0.2em}0.90 &	\hspace{-0.2em}0.90\\
\hline 
\end{tabular}\end{center}
\caption{Simulation: mean squared error (MSE), standard error of the MSE (SE), relative efficiency (REFF) and asymptotic relative efficiency (ARE) of the estimates in the uncontaminated case for $n=100$ and $p=2$.}
\label{tabla3}
\end{table}

\subsection{Results}

Table \ref{tabla3} displays the mean squared errors, the standard errors and the relative efficiencies and asymptotic relative efficiencies with respect to the MLE  for the uncontaminated case. It is seen that the relative efficiencies of all robust estimates (computed as the ratio of their respective MSEs and the MSE of the MLE) are close to their asymptotic values. The $\tau$- and MM- estimates
have similar high efficiencies, and both outperform the S-estimator.

\begin{figure}[t]
\begin{center}
\scalebox{.55}{\includegraphics{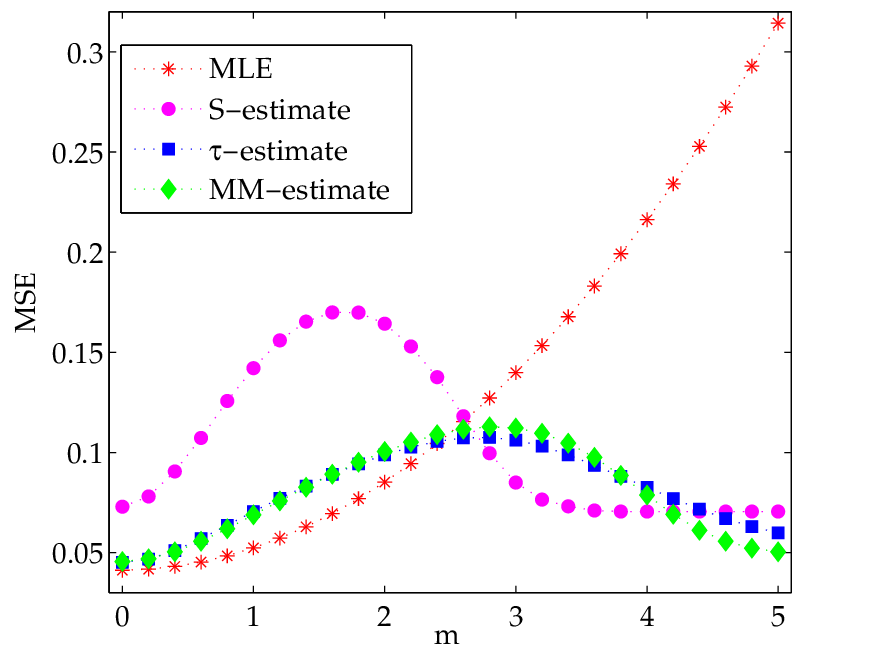}}
\end{center}
\caption{Simulation: mean squared errors for $q=2$ and $x_0=1$.}
\label{figura1}
\end{figure}

\begin{figure}[h!]
\begin{center}
\scalebox{.55}{\includegraphics{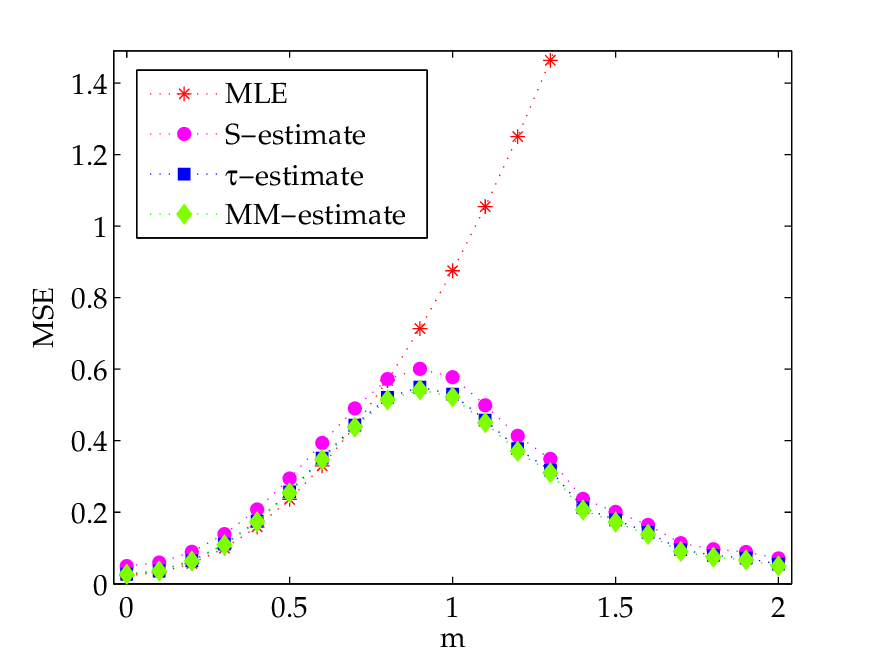}}
\end{center}
\caption{Simulation: mean squared errors for $q=2$ and $x_0=10$.}
\label{figura2}
\end{figure}

\begin{figure}[h!]
\begin{center}
\scalebox{.55}{\includegraphics{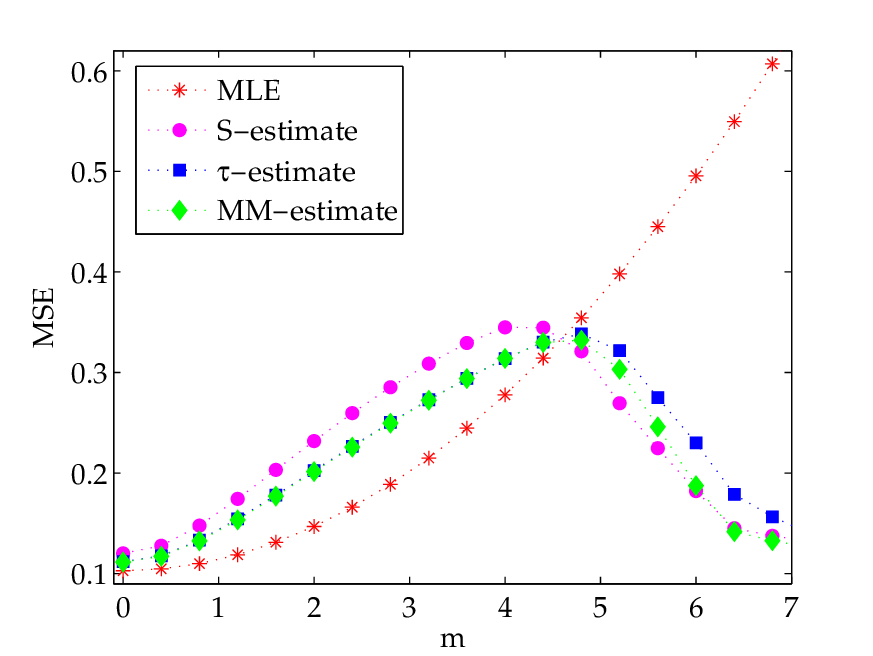}}
\end{center}
\caption{Simulation: mean squared errors for $q=5$ and $x_0=1$.}
\label{figura3}
\end{figure}
\begin{figure}[h!]
\begin{center}
\scalebox{.55}{\includegraphics{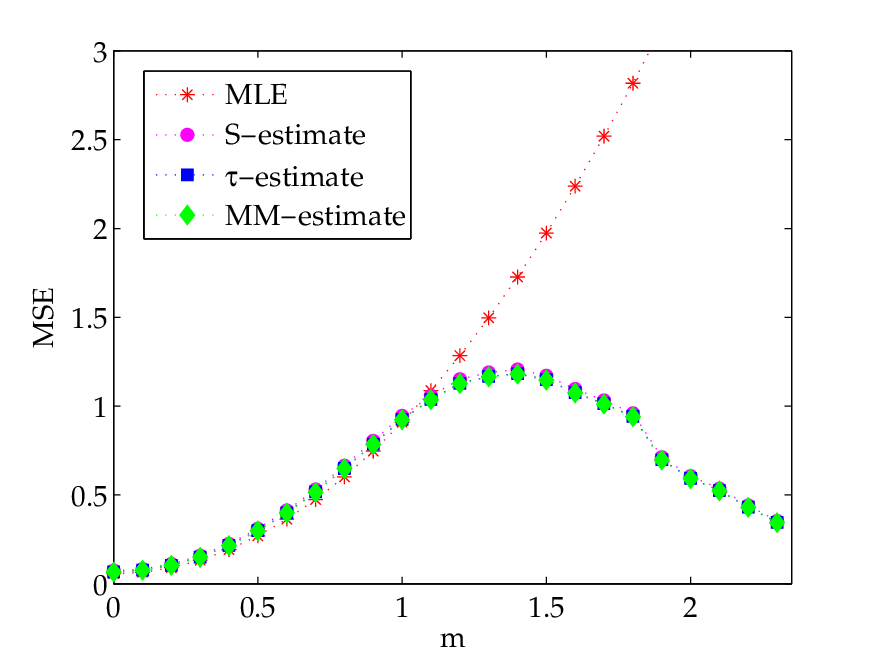}}
\end{center}
\caption{Simulation: mean squared errors for $q=5$ and $x_0=10$.}
\label{figura3b}
\end{figure}

In Figures \ref{figura1}, \ref{figura2}, \ref{figura3} and \ref{figura3b} we show the MSEs of the different estimates under contamination. 

In Figure \ref{figura1}, which corresponds to $q=2$ and $x_0=1$, we observe that the MM- and $\tau$-estimates behave similarly, both having a smaller MSE than the S-estimate except when $m$ is (approximately) between $2.8$ and $4$. In this case, the S-estimate has the largest maximum MSE among the robust estimates. As expected, the MSE of the MLE increases without bound for large $m$. Figure \ref{figura2} shows the results for $q = 2$ and $x_0 = 10$. S-, $\tau$- and MM-estimates behave similarly. In Figure \ref{figura3}, which corresponds to $q=5$ and $x_0=1$, the three robust estimates are seen to follow
essentially the same pattern. For $m\leq 4.8$ (approximately) the $\tau$- and MM-
estimates have similar behaviors, both outperforming the S-estimate. For $m >4.8$ the S- and MM-estimates have similar behaviors, both outperforming the
$\tau$-estimate. For $q = 5$ and $x_0=10$ (figure \ref{figura3b}) the behavior of the robust estimates is similar to the one observed for $q = 2$ and $x_0=10$ (figure \ref{figura2}).

\section{An example with real data}

In this Section we analyze a dataset corresponding to electron-probe X ray microanalysis of archeological glass vessels (Janssens et al., \cite{vessel}). For each of $n=180$ vessels we have a spectrum on 1920 frequencies and the contents of
13 chemical compounds; the purpose is to predict the contents on the basis of
the spectra. In order to limit the size of our data set, we considered only two
compounds (responses): P$_2$O$_5$ and PbO; and chose 12 equispaced frequencies
between 100 to 400. This interval was chosen because the values of $x_{ij}$ are almost null for frequencies below 100
and above 400. We have therefore $p=13$ and $q = 2$.

\begin{figure}[t]
\hspace{0.05cm}\scalebox{.55}{\includegraphics{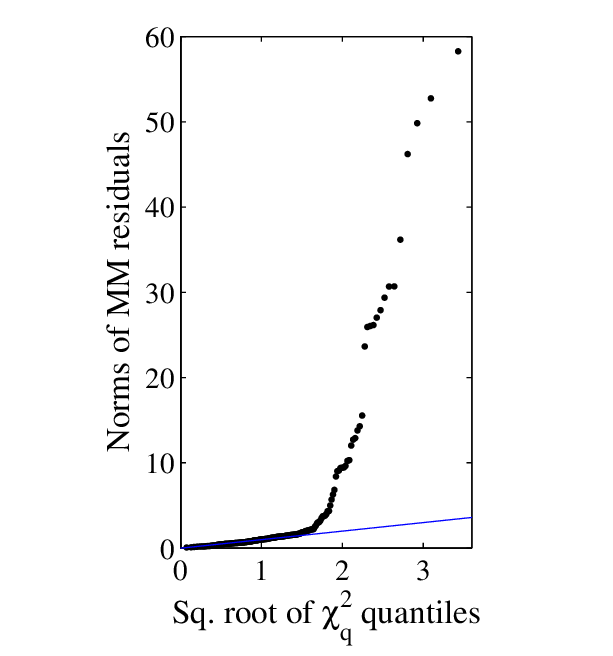}\hspace{-5em} \includegraphics{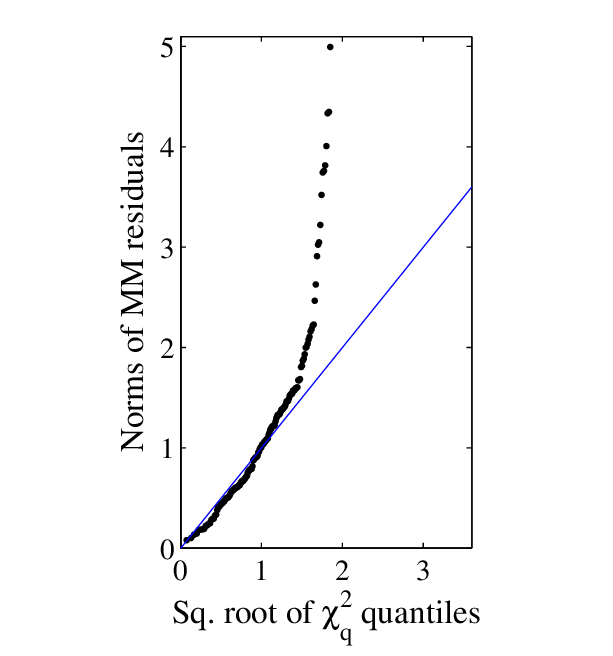}\hspace{-5em} \includegraphics{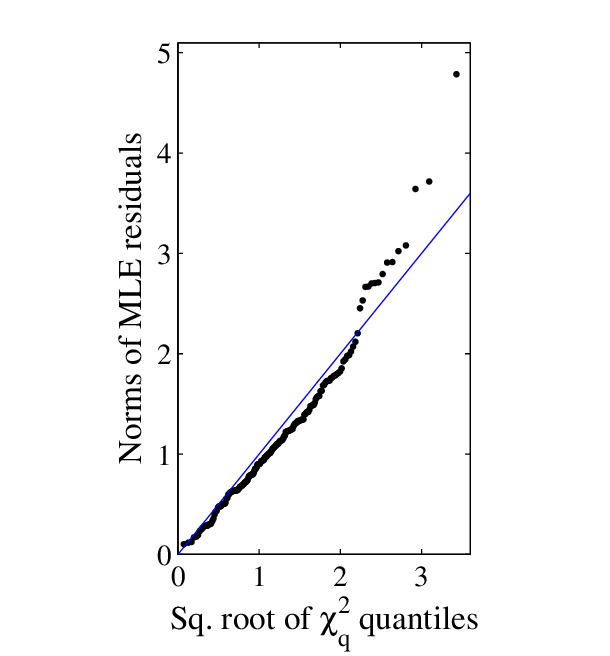}}
\caption{QQ-plots of the Mahalanobis norms of the residuals
of the MM-estimate (left), the MLE (right) and the MM-estimate in the same interval as the MLE (center).}
\label{figura4}
\end{figure}

We considered two multivariate regression estimates: the MLE and our MM-estimate. As initial estimate for the MM-estimate we use an S-estimate. The S- and the MM-estimates employ bisquare $\rho$-functions with constants such that the MM-estimate has Gaussian ARE equal to 0.95 and the S-estimate has breakdown point 0.5. In Figure \ref{figura4} we present QQ-plots of the Mahalanobis norms of the residuals
of  the MLE and the MM-estimate against the root quantiles of the chi-squared distribution with $q$ degrees of freedom. The QQ-plot of the MM-estimate shows clear outliers.

\begin{figure}[t]
\hspace{1.05cm}\scalebox{0.65}{\includegraphics{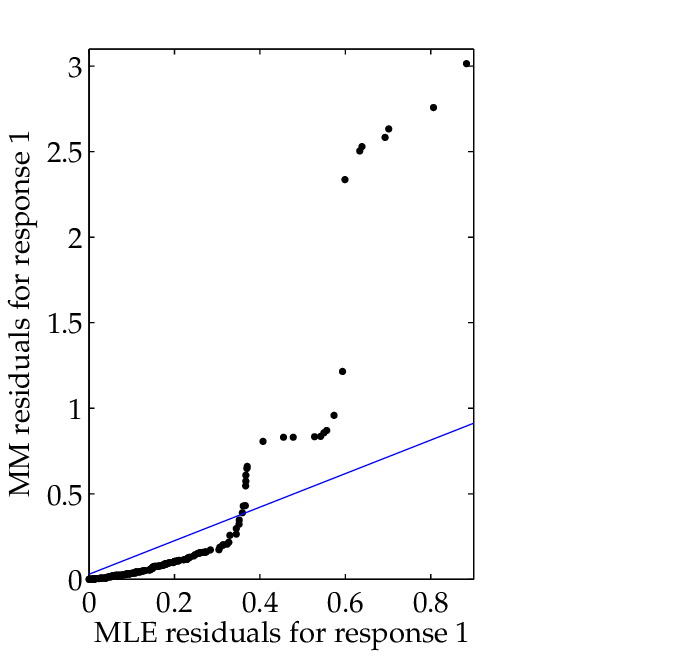}\hspace{-2em} \includegraphics{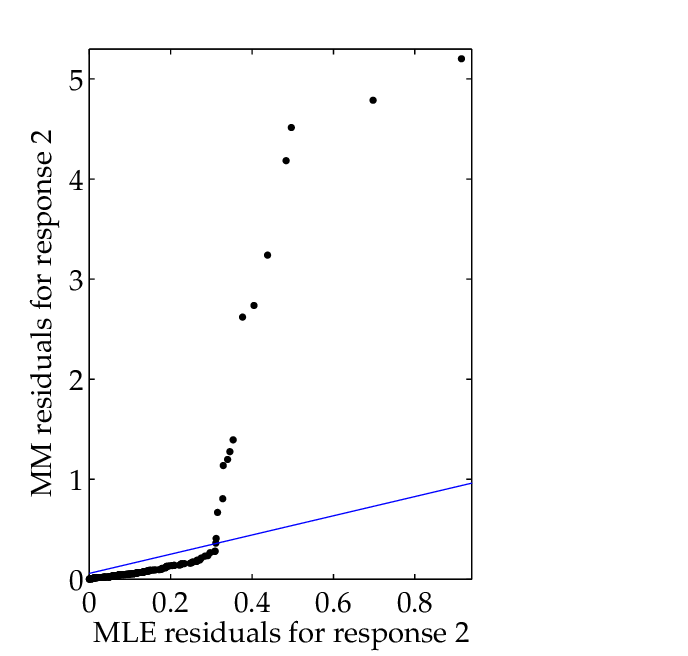}}
\caption{QQ-plots of sorted absolute residuals of MM-estimates vs sorted absolute residuals of MLE for each component of the response. Left plot corresponds to P$_2$O$_5$ (first component) and right to PbO (second component).}
\label{figure7}
\end{figure}

\begin{table}[t]
\begin{center}\begin{tabular}{ll} 
\hline
 \rowcolor[rgb]{0.8,0.8,0.8}  \hspace{1.5cm}MLE & \hspace{1cm}MM-estimate\\ 
\hline
$\left(
\begin{array}[pos]{cc}
	0.0645 & -0.0008\\
	-0.0008 & 0.0348
\end{array} \right)$ &	$\left(
\begin{array}[pos]{cc}
	0.0102 & -0.0014\\
	-0.0014 & 0.0084
\end{array} \right)$\\ 
\hline 
\end{tabular}
\end{center}
\caption{MLE and MM-estimate of the covariance matrix of the errors.}
\label{tabla4}
\end{table}

\begin{table}[t]
\begin{center}\begin{tabular}{l|ll|ll|ll|l} 
\hline
 \rowcolor[rgb]{0.8,0.8,0.8}  Criterion  & \hspace{0.5cm}MLE&\hspace{-0.5cm} & \hspace{0.1cm}$\tau$-estimate&\hspace{-0.5cm} & \hspace{0.cm}MM-estimate&\hspace{-0.5cm} & \hspace{-0.cm}MM-univariate\hspace{-1.5cm}\\ 
\hline
 \rowcolor[rgb]{0.8,0.8,0.8}  Component  & 1 &  \hspace{-0.5cm}2 & \hspace{0.cm}1 &  \hspace{-1cm}2& \hspace{0cm}1 &  \hspace{-1.3cm}2& \hspace{0cm}1\hspace{1cm}2\hspace{-1.5cm}\\ 
\hline 
MSE \hspace{0.5cm}& 0.081 &	\hspace{-0.5cm}0.051&\hspace{0.cm}0.351 &	\hspace{-1cm}0.806 &\hspace{0cm}0.340 &	\hspace{-1.3cm}0.682&\hspace{0cm}0.354\hspace{0.4cm}0.762\hspace{-1.5cm}\\ 
\hline 
$\tau$-scale \hspace{0.5cm}& 0.044 &	\hspace{-0.5cm}0.022 &\hspace{0.cm}0.005 &	\hspace{-1cm}0.007&\hspace{0cm}0.008 &	\hspace{-1.3cm}0.006&\hspace{-0cm}0.005\hspace{0.4cm}0.006\hspace{-1.5cm}\\
\hline 
\end{tabular}
\end{center}
\caption{Mean square error (MSE) and $\tau$-scale of the prediction errors of the MLE, multivariate MM-estimate, $\tau$-estimate and univariate MM-estimate for each component of the response, computed by  cross-validation.}
\label{tabla5}
\end{table}

In Figure \ref{figure7} are compared the sorted absolute values of the residuals of the MLE with those corresponding to the MM-estimator for each component of the response. 

The right and left panels of Figure \ref{figura4} show respectively the QQ-plots of the Mahalanobis norms of the residuals of the MLE and the MM-estimate against the square root quantiles of the chi-squared distribution with $q$ degrees of freedom. For ease of comparison, the center panel shows the MM's QQ-plot truncated to the size of the MLE's QQ-plot. The latter shows a very good fit of the norms to the chi-squared distribution, and therefore points out no suspect points, while the MM-estimate's QQ-plot indicates some 30 possible outliers, i.e. about 16$\%$ of the data.

The MLE's norms are in general smaller than the MM-estimate's norms, but this does not mean that the former gives a better fit, since here the residuals are normalized by the respective estimated residual dispersion matrices $\widehat{\Sig}_0$. Figure \ref{figure7} compares the sorted absolute values of the (univariate) residuals of the MLE with those of the MM-estimate for each response. We can see that the majority of the residuals corresponding to the MM-estimate are smaller than those of the MLE.

To understand why the MLE's norms are in general smaller than the MM-estimate's norms, while the respective residuals are in general smaller, we show
in Table \ref{tabla4} the estimates given by the MLE and MM-estimate of the dispersion matrix of the errors. It is seen that the former is ``much larger'' than the latter, in that its two diagonal elements are respectively six and four times those of the latter.

\begin{figure}[t]
\hspace{0.8cm}\scalebox{0.6}{\includegraphics{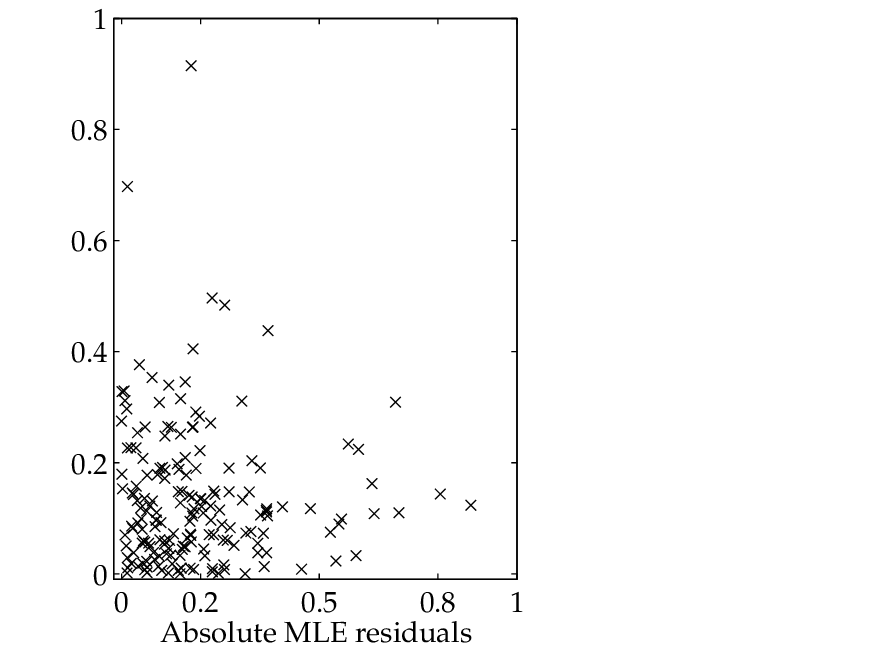}\hspace{-4cm}\includegraphics{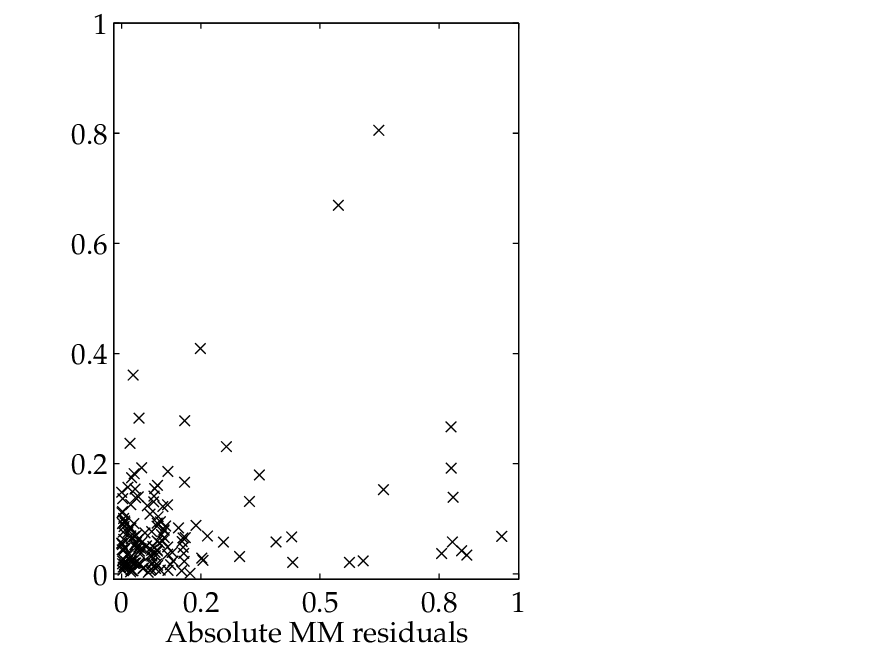}}
\caption{Absolute values of the coordinates of the bidimensional residual vectors corresponding to the MLE (left) and to the MM-estimate (right).}
\label{figure8}
\end{figure}

To complete the description of the estimates' fit, Figure \ref{figure8} shows the absolute values of the coordinates of the bidimensional residual vectors, in the right panel (corresponding to residuals of the MM-estimate) is truncated to the size of the left panel (corresponding to residuals of  the MLE), and consequently a 10$\%$ of the absolute residuals of the MM-estimate is not shown. It is seen that, while the residuals from MM have a larger range than those from the MLE, they are in general more concentrated near the origin. In general, we may conclude that the MM-estimate gives a good fit to the bulk of the data, at the expense of misfitting a reduced proportion of atypical points, while the MLE tries to fit all data points, including the atypical ones, with the cost of a poor fit to the bulk of the data.

We compared the predictive behaviors of the MLE and the MM-estimates through five-fold cross validation. We also included the univariate MM-estimates corresponding to each component of the response and the $\tau$-estimate proposed by Garc\'ia Ben et al. \cite{Tau multi}. As initial estimate for the univariate MM-estimates we use S-estimates. The $\tau$-, S- and the univariate MM-estimates employ bisquare $\rho$-functions with constants such that the univariate MM- and $\tau$-estimates have Gaussian ARE equal to 0.95 and the S-estimate has breakdown point 0.5. We considered two evaluation criteria: the mean squared error (MSE) and a robust criterion, namely a $\tau$-scale \eqref{tauscale} of the predictive errors, both computed separately for each component of the response. In the $\tau$-scale $s$ is an M-scale with breakdown point 0.5 and $\rho_2$ is a bisquare $\rho$-function with constant such that the $\tau$-scale has Gaussian asymptotic efficiency equal to 0.85. 

Table \ref{tabla5} shows the results. According to the MSE, the MLE is much better than the robust estimates. However, the $\tau$-scales yield the opposite conclusion. The reason of this fact is the MSE's sensitivity to outliers. This result shows how misleading a non-robust criterion may be.  According to the $\tau$-scale, the predictive performance of our MM-estimate for the second component is slightly better than that of the $\tau$-estimate, while the opposite occurs for the first component. The results obtained with the univariate MM-estimates are similar to those of the multivariate MM.

\begin{figure}[t]
\hspace{0.5cm}\scalebox{0.6}{\includegraphics{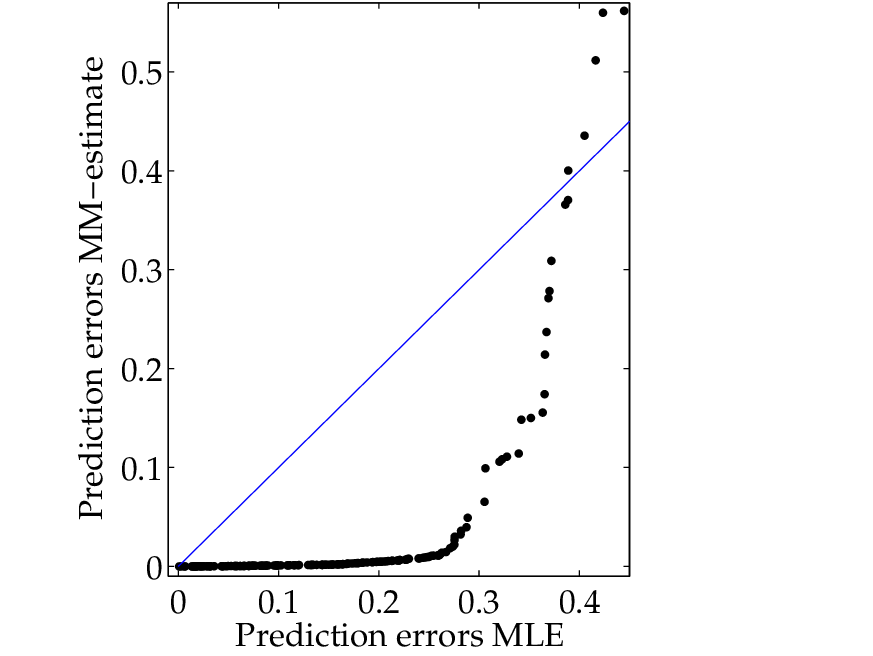}\hspace{-4cm}\includegraphics{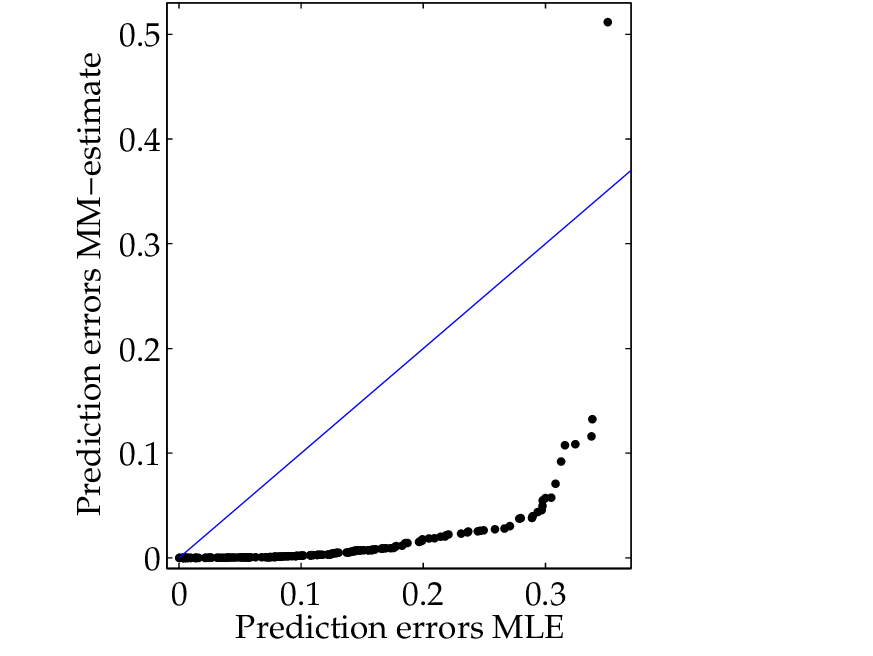}}
\caption{QQ-plots of sorted absolute prediction errors of MM-estimates vs sorted absolute prediction errors of MLE for each component of the response, computed by cross-validation. Left plot corresponds to P$_2$O$_5$ (first component) and right to PbO (second component).}
\label{fig5}
\end{figure}

The QQ-plots in Figure \ref{fig5} compare for each response component the absolute values of the sorted cross validation prediction errors of our MM-estimate with those of the MLE. For reasons of scale, in each QQ-plot the observations with the 12 largest absolute prediction errors were omitted. We can see that most points lie below the identity line representing the identity function, showing that the MM-estimate provides a better prediction for the bulk of the data.

\section{Conclusions}
In this paper we have presented MM-estimates for the multivariate linear model and showed that they maintain the same good theoretical properties
as in the univariate case, such as a high breakdown point and high Gaussian efficiency. The simulation study indicates that it has the desired high efficiency, and that its behavior is in general similar, and in several situations superior, to that of the $\tau$-estimate; it is also more efficient, and in most situations more robust, than the S-estimate. In the example with real data, our MM-estimate gives a good fit to the bulk of the data, pointing out the existence of atypical points, and shows a good predictive behavior.
\bigskip

\appendix
\section{Appendix}

Before showing some of the properties of the MM-estimate, we set the notation for norms of vectors and matrices that we will use later:

Given $\vn\in\Real^{m}$, we denote its \textit{2-norm or Euclidean norm} as:
\[
\|\vn\|=\|\vn\|_2=\left( \sum_{i=1}^{m}v_i^2\right)^{1/2},
\]
where $v_{i}$ represents the $i$th element of $\vn$.

Given a matrix $\A\in \Real^{r\times m}$, its \textit{spectral norm} $\|\cdot\|_{sp}$, is defined as follows:

\begin{align}\label{normasp}
\nonumber \|\A\| = \|\A\|_{sp}&= \max\{\|\A\vn\| : \vn\in \Real^m \mbox{ with }\|\vn\| = 1\}\\
&= \max\left\{\frac{\|\A\vn\|}{\|\vn\|} : \vn\in \Real^m \mbox{ with }\vn\ne 0\right\}.
\end{align}

Its \textit{2-norm or Frobenius norm} $\|\cdot\|_{2}$ is its Euclidean norm if we think  the matrix $\A$ as a
vector of $\Real^{m.r}$, i.e.
\begin{equation}\label{norm2}
\|\A\|_2 =\left[\text{tr}(\A'\A)\right]^{1/2}=\left(\sum_{i=1}^m\sum_{j=1}^n |a_{ij}|^2\right)^{1/2},
\end{equation}
where $a_{ij}$ represents to the $(i,j)$th element of the matrix $\A$ and $\text{tr}(\cdot)$ denotes the trace.

Given $\V\in \Real^{m\times m}$ we denote its eigenvalues as 
\[\lambda_1(\V)\geq\hdots\geq\lambda_m(\V),\]
then if $\V$ is positive definite  
\begin{equation}\label{normlambda}
\|\V\|=\max_j\{\lambda_j(\V)\}=\lambda_1(\V).
\end{equation}

\begin{rmk}\label{topo equiv}
Recall also that for any two norms $||\cdot||_{a}$ and $||\cdot||_{b}$, we have that
\[\alpha\left\|\A\right\|_a\leq\left\|\A\right\|_{b}\leq \beta\left\|\A\right\|_\alpha\]
for some $\alpha$ and $\beta$ and for all matrices $\A\in\Real^{r \times m}$. In other words, they are equivalent norms, i.e. they induce the same topology in $\Real^{r \times m}$.
For $||\cdot||_{a}=||\cdot||$ and $||\cdot||_{b}=||\cdot||_2$ we have $\alpha=1$ and $\beta=\sqrt{k}$ where $k$ is the rank of the matrix $\A$.

Recall also that spectral norm and the Frobenius norm are matrix norms, i.e., for any pair of matrices in $\Real^{m\times m}$ $\A$ and $\B$  
\begin{equation}\label{matricial}
\|\A\B\|_2\leq\|\A\|_2\|\B\|_2 \;\;\text{ and }\;\; \|\A\B\|\leq\|\A\|\|\B\|. 
\end{equation}
This property will be used in several times.
\end{rmk}

\bigskip

Before proving Theorem \ref{existencia} we will prove the following Lemma:

\begin{lem}\label{contMah}
Let $\y\in \Real^q$ and $\x\in \Real^p$ be fixed vectors.  The function
\[
d^2(\B,\Sig)=(\y-\B'\x)'\Sig^{-1}(\y-\B'\x)
\]
is continuous in $\Real^{p\times q}\times \mathcal{S}_q$.
\end{lem}

\noindent\textbf{\underline{Proof:}} We only give the main ideas of the proof. Without loss of generality, due to Remark \ref{topo equiv}, we can consider in $\Real^{p\times q}\times \Real^{q\times q}$ the topology induced by the norm 
\begin{equation}\label{trinorm}
\|\hspace{-0.1em}|(\A,\V)\|\hspace{-0.1em}|=\sup\{\|\A\|_2,\|\V\|\}.
\end{equation}

Given $(\B_0,\Sig_0)$ in $\Real^{p\times q}\times \Real^{q\times q}$, the proof consists in find an upper bound of $\left|d^2(\B,\Sig)-d^2(\B_0,\Sig_0)\right|$ that tends to $0$ when  $\|\hspace{-0.1em}|(\B,\Sig)-(\B_0,\Sig_0)\|\hspace{-0.1em}|\rightarrow 0.$

Adding and subtracting $d^2(\B_0,\Sig)$ we have that 
\begin{eqnarray*}\label{t1mast2}
\left|d^2(\B,\Sig)-d^2(\B_0,\Sig_0)\right|\hspace{-0.7em}&\leq&\hspace{-0.7em} \left|d^2(\B,\Sig)-d^2(\B_0,\Sig)\right|+\left|d^2(\B_0,\Sig)-d^2(\B_0,\Sig_0)\right|.
\end{eqnarray*}

Let $r(\B_0)=(\y-\B_0'\x)$. Using basic tools from linear algebra we obtain
\[\left|d^2(\B,\Sig)-d^2(\B_0,\Sig)\right|\leq q\|\hspace{-0.1em}|(\B,\Sig)-(\B_0,\Sig_0)\|\hspace{-0.1em}|\left(\|\x\|+2\left\|r(\B_0)\right\|\right)\|\x\|\lambda_1(\Sig^{-1}),\]
and by Weyl's Perturbation Theorem (see \cite{Bh}, pg. 63), we have that 
\begin{equation}\label{weyl}
\lambda_1(\Sig^{-1})=\frac{1}{\lambda_q(\Sig)}< \frac{1}{\lambda_q(\Sig_0)-\|\hspace{-0.1em}|(\B,\Sig)-(\B_0,\Sig_0)\|\hspace{-0.1em}|},
\end{equation}
combining these inequalities we obtain a bound of $\left|d^2(\B,\Sig)-d^2(\B_0,\Sig)\right|$ that tends to $0$ when  $\|\hspace{-0.1em}|(\B,\Sig)-(\B_0,\Sig_0)\|\hspace{-0.1em}|\rightarrow 0.$

If $r(\B_0)=\mathbf{0}$ the lemma is proved, otherwise using the Cauchy-Schwarz inequality and \eqref{weyl} we have
\begin{eqnarray*}
\left|d^2(\B_0,\Sig)-d^2(\B_0,\Sig_0)\right|\leq\frac{\|\Sig_0^{-1}\|\|r(\B_0)\|^2\|\hspace{-0.1em}|(\B,\Sig)-(\B_0,\Sig_0)\|\hspace{-0.1em}|}{\lambda_q(\Sig_0)-\|\hspace{-0.1em}|(\B,\Sig)-(\B_0,\Sig_0)\|\hspace{-0.1em}|},
\end{eqnarray*}
which completes the proof.\QED

\medskip

\noindent\textbf{\underline{Proof of Theorem \ref{existencia}:}} By Lemma \ref{contMah}, it sufices to show that there exist $t_1$ and $t_2$ such that
\begin{equation}\label{suc acotada}
\displaystyle{\inf_{(\B,\Ga)\in\mathcal{C}}\sum_{i=1}^n\rho\left(\frac{d_i(\B,\Ga)}{\sighat}\right)>\sum_{i=1}^n\rho\left(\frac{d_i(\Btil,\Gatil)}{\sighat}\right)}
\end{equation}
where 
\begin{equation}\label{C}
\mathcal{C}=\{(\B,\Ga)\in\Real^{p\times q}\times \mathcal{S}_q \text{ with } |\Ga|=1 : \lambda_q(\Ga)\leq t_1 \text{ or } \|\B\|_2\geq t_2\}.
\end{equation}

By definition of $k_n$ we have that for all $\boldsymbol{\theta}\in\Real^{p+q}$ 
\[
\sharp\{i:|\boldsymbol{\theta}'\z_i|>0\}/n\geq 1-(k_n/n).
\]
Taking 0.$5<\delta<1-(k_n/n)$ and using a compactness argument we can find $\varepsilon>0$ such that
\begin{equation}\label{infcomp}
\displaystyle{\inf_{\boldsymbol{\|\theta}\|=1}\sharp\{i:|\boldsymbol{\theta}'\z_i|>\varepsilon\}/n\geq\delta.}
\end{equation}
Let $(\B,\Ga)\in \Real^{p\times q}\times \mathcal{S}_q$ be such that $|\Ga|=1$, $\boldsymbol{\Delta}$ be the diagonal matrix of eigenvalues of $\Ga$ order from lowest to highest and $\mathbf{U}$ be the orthonormal matrix of eigenvectors of $\Ga$ which verifies $\Ga=\mathbf{U}\boldsymbol{\Delta}\mathbf{U}'$. Then
\begin{eqnarray}\label{cuenta}
\nonumber d_i^2(\B,\Ga)&=&(\y_i-\B'\x_i)'\Ga^{-1}(\y_i-\B'\x_i)\\
&=&(\mathbf{U}'\y_i-\mathbf{U}'\B'\x_i)'\mathbf{\Delta}^{-1}(\mathbf{U}'\y_i-\mathbf{U}'\B'\x_i)\geq \frac{(\mathbf{e}'\z_i)^2}{\lambda_q(\Ga)},
\end{eqnarray}
with $\mathbf{e}=(-\vn_1,\boldsymbol{\upsilon}_1)$ where $\boldsymbol{\upsilon}_1$ and $\vn_1$ are the first row of $\U'$ and $\V=\U'\B'$, respectively.

Since $\|\mathbf{e}\|\geq 1$, by \eqref{infcomp} we have at least $n\delta$ values of $d_i(\B,\Ga)$ greater or equal than $\varepsilon/\sqrt{\lambda_q(\Ga)}.$ Hence
\[
\sum_{i=1}^n\rho_1\left(\frac{d_i(\B,\Ga)}{\sighat}\right)\geq n\delta\rho_1\left(\frac{\varepsilon}{\sighat\sqrt{\lambda_q(\Ga)}}\right).
\]
Let $t$ be such that $\displaystyle{\rho_1\left(\frac{t}{\sighat}\right)=\frac{1}{2\delta_1}}$, with $0.5<\delta_1<\delta$, and let $\displaystyle{t_1=\frac{\varepsilon^2}{t^2}}$, then if $\lambda_q(\Ga)\leq t_1$ we obtain the inequality

\begin{equation}\label{parte1}
\sum_{i=1}^n\rho_1\left(\frac{d_i(\B,\Ga)}{\sighat}\right)\geq n\delta\rho_1\left(\frac{t}{\sighat}\right)>\frac{n}{2}.
\end{equation}

If $\lambda_q(\Ga)> t_1$ and $\|\B\|_2\geq t_2$, all eigenvalues of $\Ga$ are smaller than $1/t_1^{q-1}$ and at least one column of $\B$ has a norm greater o equal than $t_2/\sqrt{q}$. By \eqref{norm2}, we have
\[
\|\V\|^2_2=\|\U'\B'\|^2_2=\text{tr}(\B\U\U'\B')=\text{tr}(\B\B')=\|\B'\|^2_2=\|\B\|^2_2,
\]
and therefore exists a $k$ such that $\|\vn_k\|\geq t_2/\sqrt{q}$ where $\vn_k$ is the $k$th row of $\V$.

Then proceeding as in \eqref{cuenta} we obtain
\[
d_i^2(\B,\Ga)\geq(\mathbf{e}_k'\z_i)^2t_1^{q-1},
\]
where $\mathbf{e}_k=(-\vn_k,\boldsymbol{\upsilon}_k)$ and $\boldsymbol{\upsilon}_k$ is the $k$th row of $\U'$.

By \eqref{infcomp}, at least $n\delta$ values of $d_i(\B,\Ga)$ are greater or equal than $\varepsilon \|\mathbf{e}_k\|/\sqrt{t_1^{q-1}}$ and $\|\mathbf{e}_k\|^2=1+\|\vn_k\|^2\geq t_2^2/q$.

Then if we take $\displaystyle{t_2=\frac{t}{\varepsilon}\sqrt{{q}/{t_1^{q-1}}}}$ we obtain
\begin{equation}\label{parte2}
\sum_{i=1}^n\rho_1\left(\frac{d_i(\B,\Ga)}{\sighat}\right)\geq n\delta\rho_1\left(\frac{t}{\sighat}\right)=\frac{n}{2}.
\end{equation}

Then by \eqref{parte1} and \eqref{parte2} 
\[
\inf_{(\B,\Ga)\in\mathcal{C}}\sum_{i=1}^n\rho_1\left(\frac{d_i(\B,\Ga)}{\sighat}\right)\geq\frac{n}{2}
\]
and by \eqref{relacion-rhos} and \eqref{M-escala}
\[
\sum_{i=1}^n\rho_1\left(\frac{d_i(\Btil,\Gatil)}{\sighat}\right)\leq\sum_{i=1}^n\rho_0\left(\frac{d_i(\Btil,\Gatil)}{\sighat}\right)=\frac{n}{2},
\]
and this proves the Theorem.\QED

\bigskip
Before proving Theorem \ref{teoMMec} we will give some results on matrix derivatives that will be used later.

Let $\mathbf{b}$ be a vector and $\V$ be a symmetric matrix,
\begin{equation}\label{derivb}
\displaystyle{\frac{\partial \mathbf{b'Vb}}{\partial\mathbf{b}'}=2\mathbf{b'V}}, 
\end{equation}
if $\V$ is nonsingular,
\begin{equation}\label{derivV1}
\displaystyle{\frac{\partial |\V|}{\partial \V}=|\V|\V^{-1}}
\end{equation}
and 
\begin{equation}\label{derivV2}
\displaystyle{\frac{\partial \mathbf{b'\V^{-1}b}}{\partial \V}=-\mathbf{\V^{-1}bb'\V^{-1}}}.
\end{equation}
For further details see Chapter 17 of \cite{Seber}.

Let $(\B,\Sig)\in\Real^{p\times q}\times\mathcal{S}_q$, using $\ve{\B'\A}=(\A\otimes\Id_q)\ve{\B'}$,
it is easy to check that
\begin{equation}\label{deriv 4}
\frac{\partial \ve{\B'\A}}{\partial \ve{\B'}'}=(\A'\otimes\Id_q).
\end{equation}
From \eqref{deriv 4} and \eqref{derivb} it follows that
\begin{equation}\label{deriv 3}
\frac{\partial d(\B,\Sig)}{\partial \ve{\B'}'}=-\frac{(\y-\B'\x)'\Sig^{-1}(\x'\otimes\Id_q)}{d(\B,\Sig)}.
\end{equation}

\noindent\textbf{\underline{Proof of Theorem \ref{teoMMec}:}} The definition of MM-estimates can be reformulated, using the function $\Ga(\Sig):=\Sig/|\Sig|^{1/q}$, in the following way:
let $(\Bhat,\mathbf{C}_n)$ be any local minimum of $S^*(\B,\Sig)=S(\B,\Ga(\Sig))$
in $\Real^{p\times q}\times \mathcal{S}_q$, which satisfies
	\begin{equation*}
S^*(\Bhat,\mathbf{C}_n)\leq S^*(\widetilde{\B}_n,\widetilde{\Sig}_n).
\end{equation*} Finally, the MM-estimate $\Sighat$ is defined as
\begin{equation}\label{MM2ref}	
\Sighat=\sighat^2\Ga(\mathbf{C}_n).
\end{equation}
 
Differentiating $S^*(\B,\Sig)$ with respect to $\B$ and $\Sig$ we get 
\begin{equation}\label{DerMMec1}
\displaystyle{\sum_{i=1}^n\frac{\partial\rho_1\left(d_i(\B,\Ga(\Sig))/{\sighat}\right)}{\partial\B}{(\Bhat,\mathbf{C}_n)}=\mathbf{0}}
\end{equation}
and
\begin{equation}\label{DerMMec2} \displaystyle{\sum_{i=1}^n\frac{\partial\rho_1\left(d_i(\B,\Ga(\Sig))/\sighat\right)}{\partial\Sig}{(\Bhat,\mathbf{C}_n)}=\mathbf{0}}.
\end{equation}
By \eqref{deriv 4}, we have that
 
\begin{eqnarray*}
\displaystyle{\frac{\partial\rho_1\left(d_i(\B,\Ga(\Sig))/\sighat\right)}{\partial(\ve{\B'})'}}\hspace{-0.7em}&=&\hspace{-0.7em}  - \frac{1}{\sighat}\displaystyle{\psi_1\left(\frac{d_i(\B,\Ga(\Sig))}{\sighat}\right) \left(\frac{\uuhat_i(\B)'\Ga(\Sig)^{-1}}{d_i(\B,\Ga(\Sig))}\right)(\x_i'\otimes \Id_q)}\\
\hspace{-0.7em}&=&\hspace{-0.7em} -W\left(\frac{d_i(\B,\Ga(\Sig))}{\sighat}\right) \frac{\uuhat_i(\B)'\Ga(\Sig)^{-1}(\x_i'\otimes \Id_q)}{\sighat^2}.
\end{eqnarray*}
Then, by \eqref{MM2ref}, 
\begin{equation}\label{MMec3}
\frac{\partial\rho_1\left(d_i(\B,\Ga(\Sig))/{\sighat}\right)}{\partial(\ve{\B'})'}(\Bhat,\mathbf{C}_n)=- W\left(d_i(\Bhat,\Sighat)\right)\uuhat_i(\Bhat)'\Sighat^{-1}(\x_i'\otimes \Id_q).
\end{equation}
Using that $\ve{\Sighat^{-1}\uuhat_i(\Bhat)\x_i'}'=\uuhat_i(\Bhat)'\Sighat^{-1}(\x_i'\otimes \Id_q)$,  \eqref{MMec3} and \eqref{DerMMec1}, we can see that \eqref{MMec1} is true. 

Differentiating $\rho_1\left(d_i(\B,\Ga)/\sighat\right)$ with respect to $\Sig$, we get
\begin{equation}\label{MMec4}
\displaystyle{\frac{\partial\rho_1\left(d_i(\B,\Ga(\Sig))/\sighat\right)}{\partial\Sig}}=\frac{1}{\sighat}\psi_1\left(\frac{d_i(\B,\Ga(\Sig))}{\sighat}\right)\frac{\partial (d_i(\B,\Sig)|\Sig|^{1/(2q)})}{\partial\Sig}.
\end{equation}
From \eqref{derivV1} and \eqref{derivV2} we have
\begin{eqnarray}\label{MMec5}
\nonumber\frac{\partial (d_i(\B,\Sig)|\Sig|^{1/(2q)})}{\partial\Sig}\hspace{-0.7em}&=&\hspace{-0.7em}\displaystyle{\frac{\partial |\Sig|^{1/(2q)}}{\partial\Sig}d_i(\B,\Sig)+\frac{\partial d_i(\B,\Sig)}{\partial\Sig}|\Sig|^{1/(2q)}}\\
\nonumber \hspace{-0.7em}&=&\hspace{-0.7em}\displaystyle{ \frac{\ |\Sig|^{1/(2q)}\Sig^{-1}}{2q}d_i(\B,\Sig)\hspace{-0.1em}-\hspace{-0.1em}\frac{\Sig^{-1}\uuhat_i(\B)\uuhat_i(\B)'\Sig^{-1}}{2d_i(\B,\Sig)}|\Sig|^{1/(2q)}}\\
  \hspace{-0.7em}&=&\hspace{-0.7em}\displaystyle{ \frac{\Sig^{-1}}{2q}\left(d_i(\B,\Ga(\Sig))-q\frac{\uuhat_i(\B)\uuhat_i(\B)'\Ga(\Sig)^{-1}}{d_i(\B,\Ga(\Sig))}\right)}.
\end{eqnarray}
Then, by \eqref{MMec4} and \eqref{MMec5}, the equation \eqref{DerMMec2} results equivalent to
\[
\displaystyle{\sum_{i=1}^n \psi_1\left(\frac{d_i(\Bhat,\Ga(\mathbf{C}_n))}{\sighat}\right)\left(\frac{d_i(\Bhat,\Ga(\mathbf{C}_n))}{\sighat}-q\frac{\uuhat_i(\Bhat)\uuhat_i(\Bhat)'\Ga(\mathbf{C}_n)^{-1}}{d_i(\Bhat,\Ga(\mathbf{C}_n))\sighat}\right)}=\mathbf{0}.
\]
Rearranging and using that $W(u)=\psi(u)/u$ and that $\Sighat=\sighat^2\Ga(\mathbf{C}_n)$ we have that 
\[
\displaystyle{\Sighat\hspace{-0.2em}\left(\sum_{i=1}^n \psi_1\hspace{-0.2em}\left({d_i(\Bhat,\Sighat)}\hspace{-0.2em}\right)\hspace{-0.2em}d_i(\Bhat,\Sighat)\right)\hspace{-0.2em}=q\sum_{i=1}^n W\hspace{-0.2em}\left(\hspace{-0.2em}{d_i(\Bhat,\Sighat)}\hspace{-0.2em}\right){\uuhat_i(\Bhat)\uuhat_i(\Bhat)'}}
\]
and solving for $\Sighat$ we get \eqref{MMec2}.\QED

\bigskip

Before showing Theorem \ref{TeoBP} we will prove the following lemma:

\begin{lem}\label{lemaBP} Let $\Z=\{\z_1,\ldots, \z_n\}$, with $\z_i=(\y_i',\x_i')$ that satisfy (\ref{MLM}) and consider a $\rho$-funtion $\rho_0$.
Then the explosion breakdown point of the M-estimate of scale of the Mahalanobis norms \eqref{norma-resid}, $\sighat:=s(\dd(\Btil,\Sigtil))$, is bounded below by
\begin{equation*}\label{sigmaBP}	
\min(\varepsilon^*_n(\Z,\Btil,\Sigtil),0\text{.}5).
\end{equation*}
\end{lem}

\noindent\textbf{\underline{Proof:}} Let
\[\displaystyle{m<\min(n\varepsilon^*_n(\Z,\Btil,\Sigtil),n/2)}\]
and  let
\[\displaystyle{(\Btil^*,\Sigtil^*)=(\Btil^*(\Z^*),\Sigtil^*(\Z^*))}\] 
be an initial estimate of $(\B_0,\Sig_0)$ computed with the sample $\Z^*\in \mathcal{Z}_m$. To prove the Lemma it suffices to show that $s(\dd(\Btil^*,\Sigtil^*))$ is bounded for all $\Z^*\in \mathcal{Z}_m$.

Since $m<n\varepsilon^*_n(\Z,\Btil,\Sigtil)$, there is a compact set $K$ such that 
\[
(\Btil^*,\Sigtil^*)\in K \hspace{3em}\text{for all }\Z^*\in \mathcal{Z}_m.
\]

Then, by Lemma \ref{contMah}, there is a $t$ such that
\begin{equation}\label{sigsup}
\displaystyle{\sup_{\{i:\;\z_i=\z_i^*\}}d_i(\Btil^*,\Sigtil^*)}\leq t \hspace{3em}\text{for all }\Z^*\in \mathcal{Z}_m.
\end{equation}
Since $m/n<0.5$ we can find a $\gamma>0$ such that $m/n+\gamma<0.5$. Let $\delta$ be the value that verifies $\rho_0(\delta)= \gamma$ and let $t_0=t/\delta$.
Then using \eqref{sigsup} we have that
\begin{eqnarray*}
\displaystyle{\frac{1}{n}\sum_{i=1}^n\rho_0\left(\frac{d_i(\Btil^*,\Sigtil^*)}{t_0}\right)}&=&\displaystyle{\frac{1}{n}\sum_{\{i:\;\z_i=\z_i^*\}}\rho_0\left(\frac{d_i(\Btil^*,\Sigtil^*)}{t_0}\right)}\\
&+&\displaystyle{\frac{1}{n}\sum_{\{i:\;\z_i\neq\z_i^*\}}\rho_0\hspace{-0.2em}\left(\frac{d_i(\Btil^*,\Sigtil^*)}{t_0}\right)}\\
&\leq&\frac{(n-m)}{n}\rho_0\left({t}/{t_0}\right)+\frac{m}{n}\leq \rho_0(\delta)+\frac{m}{n}\\
&=&\gamma+ \frac{m}{n}<0\text{.}5.
\end{eqnarray*}
thus $s(\dd(\Btil^*,\Sigtil^*))\leq t_0\text{ for all }\Z^*\in \mathcal{Z}_m$ and the lemma is proved. \QED
\medskip

\noindent\textbf{\underline{Proof of Theorem \ref{TeoBP}:}} Let $\displaystyle{\varepsilon^*_n(\Z,\Btil,\Sigtil)}$ be the breakdown point of the initial estimate $(\Btil,\Sigtil)$ and  
\[\displaystyle{m<\min(n\varepsilon^*_n(\Z,\Btil,\Sigtil),[n/2]-k_n)}.\]
Let   
\[(\Bhat^*,\Sighat^*)=(\Bhat^*(\Z^*),\Sighat^*(\Z^*)) \;\;\;\; \text{ and }\;\;\;\; \displaystyle{(\Btil^*,\Sigtil^*)=(\Btil^*(\Z^*),\Sigtil^*(\Z^*))}\] 
be respectively an MM-estimate for the MLM and its initial estimate computed with the sample $\Z^*\in \mathcal{Z}_m$.

Then by (\ref{MM3}), (\ref{relacion-rhos}) and (\ref{M-escala})
\begin{eqnarray*}
\frac{1}{n}\displaystyle{\sum_{\z_i\in\Z^*}\rho_1\left({d_i(\Bhat^*,\Sighat^*)}\right)}&\leq&\displaystyle{\frac{1}{n}\sum_{\z_i\in\Z^*}\rho_1\left(\frac{d_i(\Btil^*,\Sigtil^*)}{\sighat^*}\right)}\\
&\leq&\displaystyle{\frac{1}{n}\sum_{\z_i\in\Z^*}\rho_0\left(\frac{d_i(\Btil^*,\Sigtil^*)}{\sighat^*}\right)=0\text{.}5}.
\end{eqnarray*}
Moreover, since $\sup\rho_1(u)=1$, we get 
\[
\displaystyle{\sum_{\z_i\in\Z^*}\rho_1\left({d_i(\Bhat^*,\Sighat^*)}\right)\leq \frac{n}{2}\rho_1(\infty)}.
\]

Then there exists $c<\infty$, that does not depend on $\Z^*$, such that, for at least $[n/2]$ observations of $\Z^*$, $\displaystyle{d_i^2(\Bhat^*,\Sighat^*)}<c$.

Now, since $m<[n/2]-k_n$, at least $k_n+1$ of these observations are in $\Z$, and not in a hyperplane. Then the smallest eigenvalue of $\Sighat^*$, $\lambda_q(\Sighat^*)$, is bounded below with a positive bound (for every $\x_i\in \Real^p$, the axis of the ellipsoid 
\[\{\y:(\y-\Bhat^*{}' \x_i)'\Sighat^{*-1}(\y-\Bhat^*{}'\x_i)\leq c\}\]
have lengths $\sqrt{c\lambda_j(\Sighat^{*})}; j=1,\ldots,q$. Then $\lambda_q(\Sighat^{*})>\alpha$, where $\alpha$ is a positive value not depending on $\Z^*$).

Moreover, since $|\Sighat^*|=(\sighat^*)^{2q}=s(\dd(\Btil^*,\Sigtil^*))^{2q}$, by Lemma \ref{lemaBP} the largest eigenvalue of $\Sighat^*$ is bounded above. 

To see that $\|\Bhat^*\|$ is bounded consider the set
\[
\mathcal{C}(\Bhat^*,\Sighat^*)=\{(\mathbf{v},\mathbf{w}):(\mathbf{w}-\Bhat^*{}' \mathbf{v})'\Sighat^{*-1}(\mathbf{w}-\Bhat^*{}'\mathbf{v})\leq c\}
\]
that, as we saw, contains $k_n+1$ observations of $\Z$ that are not lying on a hyperplane.

Since for symmetric matrices $\mathbf{A}$ of dimension $q\times q$, we have that
\[\lambda_q(\mathbf{A})=\inf_{\mathbf{v}} \frac{\mathbf{v}'\mathbf{A}\mathbf{v}}{\mathbf{v}'\mathbf{v}}\]
and $\lambda_q(\mathbf{A^{-1}})=1/\lambda_1(\mathbf{A})$ it follows that
\[
\|\mathbf{w}-\Bhat^*{}'\mathbf{v}\|^2\leq (\mathbf{w}-\Bhat^*{}' \mathbf{v})\Sighat^{*-1}(\mathbf{w}-\Bhat^*{}' \mathbf{v})'\lambda_1(\Sighat^*)\leq\lambda_1(\Sighat^*)c,
\]
in particular for $\mathbf{w}=\mathbf{0}$

\[
\|\Bhat^*{}'\mathbf{v}\|^2\leq \lambda_1(\Sighat^*)c.
\]
Since $\mathcal{C}(\Bhat^*,\Sighat^*)$ contains  $k_n+1$ points, there exists a constant $g$ not depending on $\Bhat^*$ or $\Sighat^*$ such that $\|\mathbf{v}\|\leq g$ implies $(\mathbf{v},\mathbf{0})\in \mathcal{C}(\Bhat^*,\Sighat^*)$. Then we have that

\[\sup_{\|\mathbf{v}\|= g}\|\Bhat^*{}'\mathbf{v}\|^2\leq \lambda_1(\Sighat^*)c,\]
that implies 

\[\|\Bhat^*{}'\|^2\leq \lambda_1(\Sighat^*)c/g^2\]
where $\|\cdot\|$ is the spectral norm defined in \eqref{normasp}.
Then, since $\|\cdot\|_2$ and $\|\cdot\|$ are equivalents, there exists a constant $\beta>0$ such that

\[
\|\Bhat^*\|_2=\|\Bhat^*{}'\|_2\leq\beta\|\Bhat^*{}'\|\leq(\beta/g)\sqrt{\lambda_1(\Gahat^*)c}
\]
for all $\Z^*\in\mathcal{Z}_m.$ This proves the Theorem. \QED\\

\bigskip

Before proving Theorem \ref{TFI} we need to prove several auxiliary Lemmas. 

\begin{lem}\label{Lema1FI}
Assume we observe $\z\in \Real^k$ with distribution $H_{\tita_1,\tita_2}$, where $\tita_1\in \Real^{m_1}$
and $\tita_2\in\Real^{m_2}$. Consider a functional M-estimate that is Fisher consistent for $\tita = (\tita_1, \tita_2)$, $\T (H) = (\T_1(H), \T_2(H))$ and an initial estimate of $\T(H)$, $\T_0 (H) = (\T_{0,1}(H), \T_{0,2}(H))$, such that 
\[E_H(h (\z, \T_1(H), \T_2(H), S(\T_0(H),H))) = 0,\]
where $h: \Real^{k+m_1+m_2+1}\rightarrow \Real^{m_1}$ is a differentiable function and  $S: \Real^{m_1+m_2}\times\mathcal{Q}\rightarrow \Real$, where $\mathcal{Q}$ is the space of distributions on $\Real^{m_1+m_2}$. Suppose that $\T$ satisfy the following strong Fisher consistency condition:
\begin{equation}\label{LFI1}
E_{H_{\tita_1,\tita_2}}(h (\z, \tita_1, \tita_2, S))=\mathbf{0} \text{ for all $S$, }
\end{equation}
and
\begin{equation}\label{LFI2}
E_{H_{\tita_1,\tita_2}}(h_3 (\z, \tita_1, \tita_2, S))=\mathbf{0} \text{ for all $S$, }
\end{equation}
where $h_i$, $1\leq i\leq 4$, is the derivative of $h$ with respect to the $i$th argument. Assume that the partial derivatives of $E_{H_{\tita_1,\tita_2}}(h_3 (\z, \tita_1, \tita_2, S(\tita,H_{\tita_1,\tita_2})))$ can be obtained differentiating with respect to each parameter inside the expectation.  Then the influence function of  $\T_1$ is given by
\begin{eqnarray*}
IF(\z_0,\T_1,\tita_1,\tita_2)&=&-\left(E_{H_{\tita_1,\tita_2}}\left(h_2(\z, \tita_1, \tita_2, S(\tita,H_{\tita_1,\tita_2}))\right)\right)^{-1}\\&&\times (h(\z_0, \tita_1, \tita_2, S(\tita,H_{\tita_1,\tita_2}))).
\end{eqnarray*}
\end{lem}

\noindent\textbf{\underline{Proof:}} Let $H_{\varepsilon}=(1-\varepsilon)H_{\tita_1,\tita_2}+\varepsilon\delta_{\z_0}$. Then $T(H_{\varepsilon})$ satisfy
\[\hspace{-5em}(1-\varepsilon)E_{H_{\tita_1,\tita_2}}(h (\z, \T_1(H_{\varepsilon}), \T_2(H_{\varepsilon}), S(\T_0(H_{\varepsilon}),H_{\varepsilon})))\] \[\hspace{5em}+ h (\z_0, \T_1(H_{\varepsilon}), \T_2(H_{\varepsilon}), S(\T_0(H_{\varepsilon}),H_{\varepsilon}))= \mathbf{0}.\]
The proof of the Lemma follows immediately differentiating the above expression with respect to $\varepsilon$ in $\varepsilon=0$ and using \eqref{LFI1} and \eqref{LFI2}. \QED\\

The following proves for the case $\Sig_0=\Id$ that the functional MM-estimates  $\T_1$ and $\T_2/S^2$ are Fisher consistent for $\B_0$ and $\Id$, respectively.

\begin{lem}\label{A.10tau}
Let $\z=(\y',\x')'$ be a random vector that satisfy the MLM (1.1) with parameters $\B_0$ and $\Sig_0=\Id$, where $\x$ satisfies (A2) and  the distribution of $\uu=\y-\B_0\x$ satisfies (A3). Let $\rho_1$ be a $\rho$-function that satisfies (A1) and $\Ga\neq\Id$ such that $|\Ga|=1$. Then 
\begin{equation*}
E_{H_0}\rho_1\left({\frac{((\y-\B'\x)'\Ga^{-1}(\y-\B'\x))^{1/2}}{\sigma}}\right)>E_{H_0}\rho_1\left(\frac{((\y-\B_0'\x)'(\y-\B_0'\x))^{1/2}}{\sigma}\right).
\end{equation*}
\end{lem}

This lemma follows immediately from Lemma A.10 of \cite{Tau multi}.
\bigskip

\begin{lem}\label{Lema2FI}
Consider the same assumptions of Theorem \ref{TFI} and suppose that $\Sig_0=\Ga_0=\Id$. Then, if $H_0$ is the distribution of $(\y',\x')'$, we have that
\[\hspace{-6em}(i)\hspace{2em} E_{H_0}\left(\frac{\partial W(d(\B_0,\Ga_0)/S)\ve{\uuhat(\B_0)\x'}}{\partial\ve{\Ga}\prime}\right)=\mathbf{0} \text{ for all } S,\]
\[\hspace{-7em}(ii)\hspace{2em} E_{H_0}\left({ W(d(\B_0,\Ga_0)/S)\ve{\uuhat(\B_0)\x'}}\right)=\mathbf{0} \text{ for all } S.\]
\end{lem}

\noindent\textbf{\underline{Proof:}} (i) By \eqref{derivV2} we have
\[ E_{H_0}\hspace{-0.2em}\left(\frac{\partial W(d(\B_0,\Ga_0)/S)\ve{\uuhat(\B_0)\x'}}{\partial\ve{\Ga}\prime}\right)\hspace{-0.2em}=\hspace{-0.2em}-E_{H_0}\hspace{-0.2em}\left(\frac{W'(\|\uu\|/S)}{2S\|\uu\|}(\x\uu'\otimes\Id)(\Id\otimes\uu\uu')\right)\]
Since the distribution of $\uu$ is assumed elliptical with $\Sig_0=\Id$, for any function $h$ we have, $E_{H_0}(h(\|\uu\|)u_iu_ju_lx_k)=0$. Then, since all the elements of the right side of the above equation have this form, part (i) of the lema is proved. (ii) follows from $E_{H_0}(h(\|\uu\|)u_ix_j)=0$ for all $i$ and $j$.\QED

\bigskip

\noindent\textbf{\underline{Proof Theorem \ref{TFI}:}} Assume $\z=(\y',\x')'$ satisfying the MLM \eqref{MLM}. Consider first the case with $\Sig_0=\Id_q$. Using Lemma \ref{Lema1FI} with $\tita_1=\ve{\B_0'}$, $\tita_2=\ve{\Ga_0}$, $S(\T_0(H),H)=S(H)$ and Lemmas \ref{A.10tau} and \ref{Lema2FI} we obtained 
\begin{eqnarray}\label{TFI1}
\nonumber IF(\z_0,\ve{\T_1'},\ve{\B_0'},\ve{\Id_q})&=&-\left(\frac{\partial E_{H_0}W\left(d(\B_0,\Ga_0)/\sigma\right)\ve{(\y-\B_0'\x)\x'}}{\partial \ve{\B'}'}\right)^{-1}\\
&&\times W\left(\|\y_0-\B_0'\x_0\|/\sigma\right)\ve{(\y_0-\B_0'\x_0)\x_0'}.
\end{eqnarray}
By \eqref{deriv 4} and \eqref{deriv 3} and the equality $\ve{\uu\x'}=(\x\otimes\Id_q)\uu$ we have that
\begin{eqnarray*}
& &\hspace{-5.5em}\frac{\partial W\left(d(\B_0,\Ga_0)/\sigma\right)\ve{(\y-\B_0'\x)\x'}}{\partial \ve{\B'}'}\\\hspace{-0.5em}&=&\hspace{-0.5em}-\frac{W'\left(d(\B_0,\sigma^2\Id_q)\right)}{d(\B_0,\sigma^2\Id_q)\sigma_0^{2}}\ve{(\y-\B_0'\x)\x'}(\y-\B_0'\x)'(\x'\otimes\Id_q)\\\hspace{-0.5em}&-&\hspace{-0.5em}W\left(d(\B_0,\sigma^2\Id_q)\right)(\x\x'\otimes\Id_q)\\
\hspace{-0.5em}&=&\hspace{-0.5em}-W'\left(\frac{\|\uu\|}{\sigma}\right)\frac{(\x\otimes\Id_q)\uu\uu'(\x'\otimes\Id_q)}{{\sigma}\|\uu\|}-W\left(\frac{\|\uu\|}{\sigma}\right)(\x\x'\otimes\Id_q).
\end{eqnarray*}

Since the distribution of $\uu$ is assumed elliptical with $\Sig_0=\Id$, for any function $h$, $E_{F_0}(h(\|\uu\|)u_iu_j)=0$ if $i\neq j$ and $E_{F_0}(h(\|\uu\|)u_i^2)=E_{F_0}(h(\|\uu\|)\|\uu\|^2)/q$. Then
\begin{eqnarray}\label{TFI2}
\nonumber & &\hspace{-3.5em}\frac{\partial E_{H_0}W\left(d(\B_0,\Ga_0)/\sigma\right)\ve{(\y-\B_0'\x)\x'}}{\partial \ve{\B'}'}\\\nonumber &=&-E_{F_0}W'\left(\frac{\|\uu\|}{\sigma}\right)\|\uu\|\frac{(E_{G_0}\x\x'\otimes\Id_q)}{q\sigma}-E_{F_0}W\left(\frac{\|\uu\|}{\sigma}\right)(E_{G_0}\x\x'\otimes\Id_q)\\&=&-\left[\frac{E_{F_0}W'\left(\|\uu\|/{\sigma_0}\right)\|\uu\|}{q\sigma}+E_{F_0}W\left(\frac{\|\uu\|}{\sigma}\right)\right](E_{G_0}\x\x'\otimes\Id_q).
\end{eqnarray}
Combining \eqref{TFI1} with \eqref{TFI2} and using then matrix equality $\ve{\mathbf{C}'\A}=(\A\otimes\Id)\ve{\mathbf{C}'}$, we obtain the proof of the Theorem in the case $\Sig_0=\Id_q$.

For the general case, let $\R$ be a matrix such that $\Sig_0=\R\R'$ and consider the following transformation $\y^*=\R^{-1}\y$. Then $\y^*={\B_0^{*}}'\x+\uu^*$, with $\uu^*=\R^{-1}\uu$ y $\B_0^*=\B_0\R'^{-1}$. Since the distribution of $\uu^*$ is given by the density \eqref{densidad} with $\Sig_0=\Id_q$ and
\[
\ve{\B_0'}=\ve{\R{\B_0^*}'},
\]
by the affine-equivariance of the estimates, we have that \[IF(\z_0,\T_1',\B_0',\Sig_0)=\R IF((\R^{-1}\y_0,\x_0),\T_1',{\B^*_0}',\Id_q), \]
where the M-scale $\sigma^*$ obtained post transforming $\y$ into $\y^*$  is $\sigma^*=\sigma|\Sig_0|^{-\frac{1}{2q}}$.\QED

\bigskip

Before proving Theorem \ref{cons} we need to prove several auxiliary Lemmas. For simplicity we will assume that the initial estimator $\Btil$ is regression- and affine-equivariant and 
$\Sigtil$ is affine-equivariant and regression-invariant. Then without loss of generality we can assume, due to Remark \ref{equiv}, that $\B_0=\mathbf{0}$ and $\Sig_0=\mathbf{I}_q$. These assumptions are not essential for the proofs.

\begin{lem}\label{cons sigma} Let $(\y_i',\x_i')$, $1\leq i\leq n$, be a random sample of the model (\ref{MLM}) with parameters $\B_0=\mathbf{0}$ and $\Sig_0=\mathbf{I}_q$, where the $\x_i$ are random and let $\rho_0$ be a $\rho$-function.  Assume that the initial estimates $\Btil$ and $\Sigtil$ are consistent for $\B_0$ and $\Sig_0$ respectively; then $\sighat$ is  consistent to $\sigma_0$ defined by the equation \eqref{4.2}.
\end{lem}

\noindent\textbf{\underline{Proof:}} Take $\varepsilon>0$, then by Lemma \ref{contMah}, we can find $\delta>0$ such that
\[
E_{H_0}\left(\inf_{\mathcal{E}}\rho_0\left(\left((\y-\B'\x)'\Sig^{-1}(\y-\B'\x)\right)^{1/2}/{(\sigma_0-\varepsilon)}\right)\right)\geq b+\delta
\]
and
\[
E_{H_0}\left(\inf_{\mathcal{E}}\rho_0\left(\left((\y-\B'\x)'\Sig^{-1}(\y-\B'\x)\right)^{1/2}/{(\sigma_0+\varepsilon)}\right)\right)\leq b-\delta
\]
where $\mathcal{E}=\{(\B,\Sig)\in\Real^{p\times q}\times \mathcal{S}_q:\|\B\|\leq \delta,\;\|\Sig - \Id_q\|\leq \delta\}.$
By the law of large numbers we have
\[
\lim_{n\longrightarrow \infty} \frac{1}{n} \sum_{i=1}^{n} \inf_{\mathcal{E}}\rho_0\left(\frac{d_i(\B,\Sig)}{(\sigma_0-\varepsilon)}\right)\geq b+\delta \;\;\;\;\;\;\;\;\;\;\;\; \text{ a.s.}
\]
and
\[
\lim_{n\longrightarrow \infty} \frac{1}{n} \sum_{i=1}^{n} \inf_{\mathcal{E}}\rho_0\left(\frac{d_i(\B,\Sig)}{(\sigma_0+\varepsilon)}\right)\leq b-\delta \;\;\;\;\;\;\;\;\;\;\;\; \text{ a.s..}
\]
Then, since  $\lim_{n\longrightarrow \infty}(\Btil,\Sigtil)=(\mathbf{0},\Id_q)$ a.s., we have
\[
\lim_{n\longrightarrow \infty} \frac{1}{n} \sum_{i=1}^{n} \rho_0\left(\frac{d_i(\Btil,\Sigtil)}{(\sigma_0-\varepsilon)}\right)\geq b+\delta \;\;\;\;\;\;\;\;\;\;\;\; \text{ a.s.}
\]
and
\[
\lim_{n\longrightarrow \infty} \frac{1}{n} \sum_{i=1}^{n} \rho_0\left(\frac{d_i(\Btil,\Sigtil)}{(\sigma_0+\varepsilon)}\right)\leq b-\delta \;\;\;\;\;\;\;\;\;\;\;\; \text{ a.s..}
\]
Therefore by the monotonicity of $\rho_0$, with probability 1 there exists $n_0$ such that for all $n\geq n_0$ we have $\sigma_0-\eps\leq\sighat\leq\sigma_0+\eps$, i.e. $\lim_{n\longrightarrow \infty}\sighat=\sigma_0$ a.s..\QED\\

The following lemma ensures the existence of a constant independent of $\B$ and $\Sig$ such that the ratio between the probability of the ellipsoid $\{(\y',\x'):(\y-\B' \x)'\Sig(\y-\B'\x)\leq \kappa\}$ and this constant is bounded by the root of each eigenvalue of $\Sig$.

\begin{lem}\label{gen L2} Suppose that the distribution of $\y$ satisfies (A3) with $\Sig_0=\Id_q$ and that $\y$ is independent of $\x$. Given $(\B,\Sig)\in \Real^{p\times q}\times \mathcal{S}_{q}$ and $\kappa>0$, consider 
\begin{equation}\label{int}
	\alpha(\B,\Sig;\kappa)=E_{H_0}I\left((\y-\B'\x)'\Sig^{-1}(\y-\B'\x)\leq\kappa\right),
\end{equation}
where $H_0$ is the distribution of $(\y',\x')$.
Then there exists a constant $\kappa_1$ independent of $\B$ and $\Sig$ such that
\[\alpha(\B,\Sig;\kappa) \leq \kappa_1 \lambda_j(\Sig)^{1/2}  \text{ for all } j, \; 1\leq j\leq q.\]
\end{lem}

\noindent\textbf{\underline{Proof:}} Note that $\V'\Sig \V=\mathbf{\Lambda}$ where $\V$ is an orthogonal matrix of $q\times q$ and $\mathbf{\Lambda}$ is a diagonal matrix whose nonzero elements are the eigenvalues of $\Sig$. Using the change of variables $\y\rightarrow \V'\y$ and (A3) we obtain that for each $j=1,\hdots, q$
\begin{eqnarray*}
E  \left(I\left(\right.\right.(\y-\B'\x)'&&\hspace{-2.5em}\Sig^{-1}(\y-\B'\x) \left.\left.\leq\kappa\right)|\x=\boldsymbol{\beta}\right)\\
&=&\int_{(\y-(\B\V)'\boldsymbol{\beta})'\mathbf{\Lambda}^{-1}(\y-(\B\V)'\boldsymbol{\beta})\leq\kappa}f_0^*(\y'\y)d\y\\ &\leq& \hspace{-0.5em} \int_{|y_j-((\B\V)'\boldsymbol{\beta})_j|\leq\sqrt{\lambda_j(\Sig)\kappa}}f_0^*\left(y_j^2+\sum_{i\neq j}y_i^2\right)\;dy_1\hdots dy_q \\ &\leq& \hspace{-0.5em}
2\sqrt{\lambda_j(\Sig)\kappa}\int f_0^*\left(\sum_{i=1}^{q-1}y_i^2\right)\;dy_1\hdots dy_{q-1}.
\end{eqnarray*}
Then if we choose
\[
\kappa_1=2\sqrt{\kappa}\int f_0^*\left(\sum_{i=1}^{q-1}y_i^2\right)\;dy_1\hdots dy_{q-1},
\]
since $\kappa_1$ does not depend on $\beta$ we obtain the desired inequality. \QED

\begin{lem}\label{acotacion} Under the assumptions of Theorem \ref{cons}, there exist positive  constants $\delta$, $L_1$ ans $L_2$  such that 
\begin{equation}\label{4.7}
\overline{\lim}_{n\longrightarrow\infty}\: \|\Bhat\|_2 \leq L_2\;\;\;\;\;\;\;\;\;\;\;\; \text{ a.s.}
\end{equation}
and
\begin{equation}\label{4.6}
\delta\leq\underline
{\lim}_{n\longrightarrow\infty}\: \|\Gahat\|\leq\overline{\lim}_{n\longrightarrow\infty}\: \|\Gahat\| \leq L_1\;\;\;\;\;\;\;\;\;\;\;\; \text{ a.s.}
\end{equation}
with $\Gahat=|\Sighat|^{-1/q}\Sighat$.
\end{lem}

\noindent\textbf{\underline{Proof:}} Let $P$ be the measure on $\Real^{q}\times \Real^{p}$ whose density is the product of $f_0(\uu)$ given in \eqref{densidad} and the density of $\x_i$, $g_0(\x)$. According to Theorem 4.2 of Ranga Rao \cite{Rao} we have
\begin{equation}\label{4.8}
\lim_{n\longrightarrow \infty} \sup_{C\subset\Real^{p+q},  \;C\text{ convex}}|P_n(C)-P(C)|=0\;\;\;\;\;\;\;\;\;\;\;\; \text{ a.s.}
\end{equation}
where $P_n$ is the empirical measure induced by the sample.

By Lemma \ref{cons sigma} there exist $n_0$ and $\sigma_1$ such that 
\begin{equation}\label{4.9}
	\sigma_1>\sighat
\end{equation}
for all $n\geq n_0$. If we consider the set 
\[
\mathcal{E}_n=\{(\y,\x):\frac{(\y-\Bhat'\x)'\Gahat^{-1}(\y-\Bhat'\x)}{\sigma_1^2}\leq\kappa\},
\]
where $\kappa$ is the constant that appears in (A1),
by \eqref{4.8} we can conclude that for large enough $n$
\[
P(\mathcal{E}_n)>P_n(\mathcal{E}_n)-b/2
\]
almost surely.

By \eqref{relacion-rhos}, \eqref{4.9} and \eqref{MM3} we have 
\begin{equation}\label{4.10}
\frac{1}{n}\sum_{i=1}^n \rho_1\left(\frac{d_i(\Bhat,\Gahat)}{\sigma_1}\right)\leq \frac{1}{n}\sum_{i=1}^n \rho_0\left(\frac{d_i(\Btil,\Sigtil)}{\sighat}\right)=b,
\end{equation}
then by (A1), 
\[\displaystyle{P_n(\mathcal{E}_n)=\frac{1}{n} \sharp{\{(\y_i,\x_i):(\y_i-\Bhat'\x_i)'\Gahat^{-1}(\y_i-\Bhat'\x_i)/\sigma_1^2\leq\kappa\}}}\geq b\]
and therefore $P(\mathcal{E}_n)> b/2$ almost surely for $n$ large enough.

By Lemma \ref{gen L2} $P(\mathcal{E}_n)\leq \lambda_j(\Gahat)^{1/2}\sigma_1\kappa_1$ for all $1\leq j\leq n$, then if $\delta=b^2/(4\kappa_1^2\sigma_1^2)$ for $n$ large enough we have that $\lambda_j(\Gahat)\geq \delta$, almost sure, for all $j$, in particular for $\lambda_1(\Gahat)=\|\Gahat\|$. Then since $|\Gahat|=1$ we have that there is a constant $L_1>0$ such that for $n$ large enough $\|\Gahat\|\leq L_1$.

By \eqref{4.10} and Lemma \ref{cons sigma} to prove \eqref{4.7} it would be enough to show that for any $\sigma>0$ there exist $L_2$ and  $\eta>0$ such that
\begin{equation}\label{4.11}
\displaystyle{\lim_{n\longrightarrow\infty} \inf_{\; \|\B\|>L_2\;}\frac{1}{n}\sum_{i=1}^n \rho_1\left(\frac{d_i(\B,\Gahat)}{\sigma}\right)\geq b+\eta} \;\;\;\;\;\;\;\; \text{a.s.}.
\end{equation}

By the Lebesgue dominated convergence Theorem, it is easy to show that for any $\sigma>0$
\begin{equation}\label{4.12}
\lim_{M_1\longrightarrow\infty}E_{F_0}\rho_1\left(\frac{\|\y\|-M_1}{L_1^{1/2}\sigma}\right)=1.
\end{equation}

By (A2), there exist $\varphi>0$, $\gamma>0$ and a finite number of sets $\mathcal{C}_1,\mathcal{C}_2,\ldots,\mathcal{C}_s$ included in $\Real^{p\times q}$ such that
\begin{equation}\label{4.13}
\bigcup_{i=1}^s \mathcal{C}_i\supset\mathcal{C}=\{\B\in\Real^{p\times q}: \|\B\|_2=1\}
\end{equation}
and
\begin{equation}\label{4.14}
P_{G_0}(\inf_{\;\B\in\mathcal{C}_i\;}\|\B'\x\|\geq \varphi)\geq b+\gamma.
\end{equation}
By \eqref{4.12} we can find $M_1$ and $\eta>0$ such that 
\begin{equation}\label{4.15}
(b+\gamma)E_{F_0}\left(\rho_1\left(\frac{\|\y\|-M_1}{L_1^{1/2}\sigma}\right)\right)>b+2\eta.
\end{equation}
Theb by \eqref{4.14} and \eqref{4.15} we have
\begin{equation}\label{4.16}
E\inf_{\;\B\in\mathcal{C}_i\;}I(\|\B'\x\|\geq \varphi)\rho_1\left(\frac{\|\y\|-M_1}{L_1^{1/2}\sigma}\right)\geq b+2\eta.
\end{equation}
Let $M_2$ be such that 
\begin{equation}\label{4.17}
P_{F_0}\left(\|\y\|\geq M_2\right)<\eta,
\end{equation}
take $M=\max\{M_1,\;M_2\}$ and $L_2=M/\varphi$, then by \eqref{4.6} and \eqref{4.13} we have

\begin{eqnarray*}
\inf_{\; \|\B\|_2>L_2\;}\hspace{-2.5em}&&\frac{1}{n}\sum_{i=1}^n \rho_1\left(\frac{d_i(\B,\Gahat)}{\sigma}\right)\geq
\inf_{\; \|\B\|_2>L_2\;}\frac{1}{n}\sum_{i=1}^n \rho_1\left(\frac{\|\y_i-\B'\x\|}{L_1^{1/2}\sigma}\right)\\
&\geq& \hspace{-0.7em}
\inf_{\; \|\B\|_2=1\;}\frac{1}{n}\sum_{i=1}^n \rho_1\left(\frac{\|\y_i\|-L_2\varphi}{L_1^{1/2}\sigma}\right)I(\|\B'\x\|>\varphi)I(\|\y_i\|<L_2\varphi)\\&\geq& \hspace{-0.7em}
\inf_{1\leq j\leq s}\inf_{\; \B\in\mathcal{C}_j\;}\frac{1}{n}\sum_{i=1}^n \rho_1\left(\frac{\|\y_i\|-M}{L_1^{1/2}\sigma}\right)I(\|\B'\x\|>\varphi)I(\|\y_i\|<M)\\&\geq&\hspace{-0.7em}
\inf_{1\leq j\leq s}\frac{1}{n}\sum_{i=1}^n\inf_{\; \B\in\mathcal{C}_j\;} \rho_1\left(\frac{\|\y_i\|-M}{L_1^{1/2}\sigma}\right)I(\|\B'\x\|>\varphi)(1-I(\|\y_i\|\geq M))
\end{eqnarray*}
Finally, using the Law of Large Numbers, \eqref{4.16} and \eqref{4.17}  we get \eqref{4.11} and this proves \eqref{4.7}.\QED

\begin{lem}\label{gen lema 4.2} 
Let $\mathbf{g}:\Real^{k}\times(\Real^{m\times n}\times\Real^{r\times t})\longrightarrow\Real$ continuous and let $Q$ be a  probability distribution on $\Real^{k}$ such that for some $\delta>0$ we have
\[
E_Q\left(\sup_{\|\hspace{-0.05em}|(\A,\V)-(\A_0,\V_0)\|\hspace{-0.05em}|\leq\delta}|\mathbf{g}\left(\z,(\A,\V)\right)|\right)<\infty,
\]
where $\|\hspace{-0.1em}|\cdot\|\hspace{-0.1em}|$ is the norm defined in \eqref{trinorm}. Let $(\Ahat_n,\Vhat_n)$ be a sequence of estimates in $\Real^{m\times n}\times\Real^{r\times t}$ such that $\lim_{n\longrightarrow\infty}(\Ahat_n,\Vhat_n)=(\A_0,\V_0)$ a.s.. Then if $\z_1,\ldots,\z_n$, are i.i.d. random variables in $\Real^{k}$ with distribution $Q$, we have
\[
\lim_{n\longrightarrow\infty}(1/n)\sum_{i=1}^n\mathbf{g}(\z_i,(\Ahat_n,\Vhat_n))=E_Q\mathbf{g}(\z,(\A_0,\V_0)) \;\;\;\;\;\;\text{a.s}.. 
\]
\end{lem}

\noindent\textbf{\underline{Proof:}} To prove the Lemma it suffices to show that for any $\varepsilon>0$ there exists $\eta>0$ such that 
\begin{equation}\label{4.18}
\overline{\lim}_{n\longrightarrow\infty} \sup _{\|\hspace{-0.05em}|(\A,\V)-(\A_0,\V_0)\|\hspace{-0.05em}|\leq\eta}(1/n)\sum_{i=1}^n\mathbf{g}(\z_i,(\A,\V)) \leq E_Q\mathbf{g}(\z,(\A_0,\V_0))+ \varepsilon
\end{equation} 
and
\begin{equation}\label{4.19}
\underline{\lim}_{n\longrightarrow\infty} \inf _{\|\hspace{-0.05em}|(\A,\V)-(\A_0,\V_0)\|\hspace{-0.05em}|\leq\eta}(1/n)\sum_{i=1}^n\mathbf{g}(\z_i,(\A,\V)) \geq E_Q\mathbf{g}(\z,(\A_0,\V_0))- \varepsilon
\end{equation}

By the Lebesgue dominated convergence Theorem we can take $0<\eta<\delta$ such that 
\[
E(\sup _{\|\hspace{-0.05em}|(\A,\V)-(\A_0,\V_0)\|\hspace{-0.05em}|\leq\eta}\mathbf{g}(\z,(\A,\V)))\leq  E_Q\mathbf{g}(\z,(\A_0,\V_0))+ \varepsilon.
\]
Then using the Law of Large Numbers we obtain 
\[
\overline{\displaystyle\lim_{n \to\infty}{}} \frac{1}{n}\sum_{i=1}^n\sup _{\|\hspace{-0.05em}|(\A,\V)-(\A_0,\V_0)\|\hspace{-0.05em}|\leq\eta}\hspace{-0.1em}\mathbf{g}(\z_i,(\A,\V)) \hspace{-0.2em}=\hspace{-0.2em}E(\sup _{\|\hspace{-0.05em}|(\A,\V)-(\A_0,\V_0)\|\hspace{-0.05em}|\leq\eta}\hspace{-0.1em}\mathbf{g}(\z,(\A,\V)))
\]
and get \eqref{4.18}. A similar procedure is performed to prove \eqref{4.19}. \QED
\medskip

\noindent\textbf{\underline{Proof of Theorem \ref{cons}:}} Consider
\begin{eqnarray*}
\mathcal{C}(\delta,L_1,L_2)\hspace{-0.5em}&=&\hspace{-0.5em}\{(\B,\Ga)\in \Real^{p\times q}\times \mathcal{S}_q: \delta\leq\|\Ga\|\leq L_1,|\Ga|=1 \text{ and } \|\B\|_2\leq L_2\},\\
\mathcal{C}_1(\varepsilon)\hspace{-0.5em}&=&\hspace{-0.5em}\{(\B,\Ga)\in \mathcal{C}(\delta,L_1,L_2):  \|\B\|_2\geq \varepsilon\} \text{ and} \\
\mathcal{C}_2(\varepsilon)\hspace{-0.5em}&=&\hspace{-0.5em}\{(\B,\Ga)\in \mathcal{C}(\delta,L_1,L_2):  \|\Ga-\Id_q\|\geq \varepsilon\}.
\end{eqnarray*}
According to the Lemmas \ref{cons sigma} and \ref{acotacion} and \eqref{MM3}, it would be enough to show that given $\varepsilon_1>0$, $\varepsilon_2>0$ and $L_1$ and $L_2$ arbitrarily large, there exist $\gamma>0$ and $\sigma_1>\sigma_0$ such that  
\begin{equation}\label{4.20}
\underline{\lim}_{n\longrightarrow\infty} \inf _{(\B,\Ga)\in \mathcal{C}_1(\eps_1)}\frac{1}{n}\sum_{i=1}^n \rho_1\left(\frac{d_i(\B,\Ga)}{\sigma_1}\right)\geq E_{F_0}\rho_1\left(\frac{(\uu'\uu)^{1/2}}{\sigma_0}\right)+\gamma \;\;\;\;\;\;\;\; \text{a.s.,}
\end{equation}
\begin{equation}\label{4.21}
\underline{\lim}_{n\longrightarrow\infty} \inf _{(\B,\Ga)\in \mathcal{C}_2(\eps_2)}\frac{1}{n}\sum_{i=1}^n \rho_1\left(\frac{d_i(\B,\Ga)}{\sigma_1}\right)\geq E_{F_0}\rho_1\left(\frac{(\uu'\uu)^{1/2}}{\sigma_0}\right)+\gamma \;\;\;\;\;\;\;\; \text{a.s.}
\end{equation}
and
\begin{equation}\label{4.22}
\lim_{n\longrightarrow\infty} \frac{1}{n}\sum_{i=1}^n \rho_1\left(\frac{d_i(\Btil,\Sigtil)}{\sighat}\right)= E_{F_0}\rho_1\left(\frac{(\uu'\uu)^{1/2}}{\sigma_0}\right) \;\;\;\;\;\;\;\;\;\; \text{a.s..}
\end{equation}
By Lemma \ref{A.10tau} we have 
\begin{equation}\label{cons fisher}
E\rho_1\left(\frac{d(\B,\Ga)}{\sigma_0}\right)>E\rho_1\left(\frac{(\uu'\uu)^{1/2}}{\sigma_0}\right)
\end{equation}
for all $\B\in \Real^{p\times q}$ and $\Ga\in \mathcal{S}_q$ with $|\Ga|=1$  such that $\Ga\neq \Id_q$.

By Lemma \ref{contMah}, \eqref{cons fisher} and the Lebesgue dominated convergence Theorem, using a standard compactness argument we can find $\sigma_1>\sigma_0$, $\gamma>0$ and a finite number of sets, $\mathcal{C}_1,\hdots,\mathcal{C}_s$, such that 
\begin{equation}\label{4.24}
E_{H_0}\inf_{(\B,\Ga)\in\mathcal{C}_j}\rho_1\left(\frac{d(\B,\Ga)}{\sigma_1}\right)>E_{F_0}\rho_1\left(\frac{(\uu'\uu)^{1/2}}{\sigma_0}\right)+\gamma
\end{equation}
and
\begin{equation}\label{4.25}
\bigcup_{j=1}^s\mathcal{C}_j\supset\mathcal{C}_1(\eps_1).
\end{equation}
By \eqref{4.25} we have 
\begin{eqnarray*}
\lim_{n\longrightarrow\infty}\inf_{(\B,\Ga)\in\mathcal{C}_1(\eps_1)}\frac{1}{n}\sum_{i=1}^n&&\hspace{-2.5em}\rho_1\left(\frac{d_i(\B,\Ga)}{\sigma_1}\right)\\&\geq& \inf_{1\leq j\leq s}\lim_{n\longrightarrow\infty}\frac{1}{n}\sum_{i=1}^n\inf_{(\B,\Ga)\in\mathcal{C}_j}\rho_1\left(\frac{d_i(\B,\Ga)}{\sigma_1}\right).
\end{eqnarray*}
Then by \eqref{4.24} and the Law of Large Numbers we get \eqref{4.20}. \eqref{4.21} is proved similarly to \eqref{4.20} and \eqref{4.22} is a consequence of Lemma \ref{gen lema 4.2}.  \QED

\bigskip

Next we will give some definitions and lemmas that will be necessary to prove the asymptotic normality of MM-estimates $\Bhat$.

\begin{Def}\label{envolvente}
Let $\mathfrak{F}$ be a class of real-valued functions on a set $\mathfrak{X}$. An \textit{envelope} for $\mathfrak{F}$ is any function $F$ such that $|f|\leq F$ for all $f$ in $\mathfrak{F}$.
\end{Def}

If $\mu$ is a measure on $\mathfrak{X}$ for which $F$ is integrable, it is natural to think of $\mathfrak{F}$ as a subset of $\mathfrak{L}^1(\mu)$, the space of all $\mu$-integrable functions. This space is equipped with a distance defined by the $\mathfrak{L}^1(\mu)$ norm. Then the closed ball with center $f_0$ and radius $R$ consists of all $f$ in $\mathfrak{L}^1(\mu)$ for which $\int|f-f_0|d\mu\leq R$.

\begin{Def}\label{euclideana}
The class $\mathfrak{F}$ is \textit{Euclidean} for the envelope $F$ if there exist positive constants $a$ and $r$ with the following property: if $0\leq \varepsilon \leq 1$ and if $\mu$ is any measure for which $\int Fd\mu<\infty$, then there are functions $f_1,\hdots, f_m$ in $\mathfrak{F}$ such that 
\begin{description}
\item[(i)] $m\leq a\varepsilon^{-r}$,
\item[(ii)] $\mathfrak{F}$ is covered by the union of the closed balls with radius $\varepsilon\int Fd\mu$ and centers $f_1,\hdots, f_m$.
\end{description} 
\end{Def}
%
In order to prove the following lemma we need Lemma 2.13 from Pakes and Pollard \cite{Pakes Pollard}. This is stated below:

\begin{lem}\label{poll2}
Let $\mathfrak{F}=\{f(\cdot,\xi):\xi\in \mathcal{C}\}$ be a class of functions on $\mathfrak{X}$ indexed by a bounded subset $\mathcal{C}$ of $\Real^d$. If there exists an $\alpha>0$ and a nonnegative function $\varphi(\cdot)$ such that
\[
|f(x,\xi)-f(x,\xi^*)|\leq\varphi(x)\|\xi-\xi^*\|^{\alpha} \;\text{ for } \;x\in \mathfrak{X} \;\;\;\text{ and }\;\;\; \xi,\xi^*\in \mathcal{C},
\]
then $\mathfrak{F}$ is Euclidean for the envelope $|f(\cdot,\xi_0)|+R\varphi(\cdot)$, where $\xi_0$ is an arbitrary point of $\mathcal{C}$ and $R=(2\sqrt{d}\sup_\mathcal{C}\|\xi-\xi_0\|)^{\alpha}$.
\end{lem}

The proof of Lemma \ref{poll2} can be found in \cite{Pakes Pollard}.

\begin{lem}\label{hipA6}
If (A4), (A5) and (A6) hold, then there exists a function $\tita(\xi)$, that to each $\xi$ in $\Real^{qp}\times\ve{\mathcal{S}_q}$ assigns a pair $(\B,\Sig)$ in $\Real^{q\times p}\times \mathcal{S}_q$, and a bounded subset $\mathcal{C}$ of $\Real^{qp}\times\ve{\mathcal{S}_q}$, such that $(\B_{0},\sigma_0^2\Sig_0)\in\tita(\mathcal{C})^{\circ}$, for which each of the classes of functions 
\begin{equation}\label{familia F}
\mathfrak{F}_{kj}=\{\phi_{kj}(\z;\tita(\xi)):\xi\in\mathcal{C}\},
\end{equation}
where $\phi_{kj}(\z;\tita)=W\left(d(\B,\Sig)\right)(y_k-\bn_k'\x)x_j$ and $\bn_k$ is the $k$th column vector of the matrix $\B$, is Euclidean for certain envelope $F_{kj}$ with $E_{H_0}F_{kj}^2<\infty$.
\end{lem}

\noindent\textbf{\underline{Proof:}} 
For each $\xi$ in $\Real^{qp}\times\ve{\mathcal{S}_q}$ there exists a unique pair $\left({\B},{\Sig}\right)$ in $\Real^{q\times p}\times \mathcal{S}_q$ such that $\xi=\left(\ve{\B}',\ve{\Sig^{-1/2}}'\right)$, then define the function $\tita(\cdot)$ as follows:
$
\tita\left(\left(\ve{\B}',\ve{\Sig^{-1/2}}'\right)\right)=\left({\B},{\Sig}\right).
$

Let $\varepsilon>0$ and $\delta=2\|\sigma_0^{-1}{\Sig_0^{-1/2}}\|$, considering the norm defined in \eqref{trinorm}, we denote by $\mathcal{B}_{\varepsilon}$ to the ball of radius $\varepsilon$ and center $({\B_{0}},\sigma_0^2{\Sig_0})$, then define
\[
\mathcal{C}=\tita^{-1}\left(\left\{({\B},{\Sig})\in\Real^{q\times p}\times \mathcal{S}_q: ({\B},{\Sig})\in\mathcal{B}_{\varepsilon} \text{ and } \|\Sig^{-1/2}\|\leq\delta\right\}\right).
\]

Let $\xi$ and $\xi^*$ be any two elements of $\mathcal{C}$ such that $\tita(\xi)=(\B,\Sig)$ and $\tita(\xi^*)=(\B^*,\Sig^*)$,  by the Mean Value Theorem there is a value $c$ between
$d(\B,\Sig)$ and $d(\B^*,\Sig^*)$ such that
\begin{equation}\label{TVM}
|W\left(d(\B^*,\Sig^*)\right)-W\left(d(\B,\Sig)\right)|=
|W'(c)||d(\B,\Sig)-d(\B^*,\Sig^*)|.
\end{equation}
Since $W$ and its derivative are continuous and with compact support there exists a constant $M$ such that $|W(u)|\leq M$ and $|W'(u)|\leq M$ for all $u$, using this and \eqref{TVM} we have
\begin{eqnarray*}
|\phi_{kj}(\z;\tita(\xi))-\phi_{kj}(\z;\tita(\xi^*))|
&\leq&|W'(c)||d(\B,\Sig)-d(\B^*,\Sig^*)||(y_k-{\bn^*_k}'\x) x_j|\\
&+&|W\left(d(\B,\Sig)\right)||{\bn^*_k}'\x x_j-\bn_k'\x x_j|\\
&\leq&M\left\{|d(\B,\Sig^*)-d(\B^*,\Sig)|\left(|y_k|+\|{\bn^*_k}\|\|\x\|\right)\right. \\
&+&\left.||{\B^*}'-\B'||_2\|\x \|\right\}|x_j|\\
&\leq&
M|x_j|\left\{|d(\B,\Sig^*)-d(\B^*,\Sig)|\left(|y_k|+\varepsilon\|\x\|\right)\right. \\
&&\left.+\|\xi-\xi^*\|\|\x \|\right\}.
\end{eqnarray*}
Applying inequalities of matrix norms we have
\begin{eqnarray*}\label{difdist}
\nonumber|d(\B,\Sig)&&\hspace{-2.2em}-\hspace{0.2em}d(\B^*,\Sig^*)|\leq\|\Sig^{-1/2}\left(\y-{\B}'\x\right)-\Sig^{*-1/2}\left(\y-{\B^*}'\x\right)\|\\
\nonumber&\leq&\hspace{-0.4em}\|\Sig^{-1/2}-\Sig^{*-1/2}\|_2\|\y\|+\|\Sig^{-1/2}{\B}'\x-\Sig^{*-1/2}{\B^*}'\x\|\\
\nonumber&\leq&\hspace{-0.4em}\|\Sig^{-1/2}-\Sig^{*-1/2}\|_2\|\y\|+\|\Sig^{-1/2}\|\|\B'-{\B^*}'\|_2\|\x\|\\
\nonumber&+&\hspace{-0.4em}\|{\Sig^{*-1/2}}-\Sig^{-1/2}\|_2\|{\B^*}\|_2\|\x\|\\
\nonumber&\leq&\hspace{-0.4em}\|\Sig^{-1/2}-\Sig^{*-1/2}\|_2\|\y\|+\delta\|\B'-{\B^*}'\|_2\|\x\|\\
\nonumber&+&\hspace{-0.4em}\|\Sig^{-1/2}-\Sig^{*-1/2}\|_2\varepsilon\|\x\|\\
&\leq&\hspace{-0.4em}\left\{(\varepsilon+\delta)\|\x\|+\|\y\|\right\}\|\hspace{-0.2em}\left(\hspace{-0.1em}\ve{\B\hspace{-0.1em}-\hspace{-0.1em}{\B^*}}',\ve{\Sig^{-1/2}-{\Sig^{*-1/2}}}'\right)\hspace{-0.2em}\|\\
&=&\hspace{-0.4em}\left\{(\varepsilon+\delta)\|\x\|+\|\y\|\right\}\|\xi-\xi^*\|.
\end{eqnarray*}

Then, if we define
\begin{equation}\label{phikj}
\varphi_{kj}(\z)=M\left\{\left((\varepsilon+\delta)\|\x\|+\|\y\|\right)(|y_kx_j|+\varepsilon\|\x\| |x_j|)+\|\x \||x_j|\right\}
\end{equation}
we have that
\begin{eqnarray*}
|\phi_{kj}(\z;\tita(\xi))-\phi_{kj}(\z;\tita(\xi^*))|\leq\varphi_{kj}(\z)\|\xi-\xi^*\|.
\end{eqnarray*}
Then we can apply Lemma \ref{poll2} and conclude that $\mathfrak{F}_{kj}$ is euclidean for the envelope
\begin{equation*}
F_{kj}(\z)=|\phi_{kj}(\z,\tita(\xi_0))|+R\varphi_{kj}(\z),
\end{equation*}
with $\xi_0\in \mathcal{C}$ such that $\tita(\xi_0)=(\B_0,\sigma_0^2\Sig_0)$ and $R=2\sqrt{q(p+q)}\sup_{\mathcal{C}}\|\xi-\xi_0\|$. The proof of $E_{H_0}F_{kj}^2<\infty$ follows immediately using that $|W(u)|\leq M$ and expanding \eqref{phikj} as a sum of products, and bounding their respective means by means of (A6). \QED

Before proving Theorem \ref{normalidad} we need to state Lemma 2.16 (page 1036) of Pakes and Pollard \cite{Pakes Pollard}.

\begin{lem}\label{pollard}
Let $\mathfrak{F}$ be a \textit{Euclidean} class with envelope $F$ such that $\int F^2dP<\infty$. For each $\eta>0$ and $\varepsilon>0$ there exists a $\delta>0$ such that
\[
\lim\sup\mathbb{P}\left\{\sup_{[\delta]}|\nu_n(f_1)-\nu_n(f_2)|>\eta\right\}<\varepsilon,
\]
where $[\delta]$ represents the set of all pairs of functions in $\mathfrak{F}$ with 
\[\int(f_1-f_2)^2dP<\delta^2\] 
and $\displaystyle{\nu_n(f)={n}^{-1/2}\sum_{i=1}^{n}\left[f(\zeta_i)-\int fdP\right]}$, where $\zeta_1,\zeta_2,\dots,\zeta_n$ are independent observations sampled from the distribution $P$.
\end{lem}

The proof of Lemma \ref{pollard} can be found in \cite{Pakes Pollard}.\\

\noindent\textbf{\underline{Proof of Theorem \ref{normalidad}:}} We denote $\boldsymbol{\theta}_n=(\Bhat,\Sighat)$ and $\boldsymbol{\theta}_0=(\B_0,\sigma_0^2\Sig_0)$.

Since we assumed that the distribution of errors $\uu$ is elliptical with density of the form \eqref{densidad}, for any function $h$ we have $E_{H_0}\left(x_ju_ih(\uu'\Sig_0^{-1}\uu)\right)=0$. This implies that
\begin{equation}
E_{H_0} \frac{1}{\sigma_0^2}W\left(d(\B,\sigma_0^2\Sig_0)\right)(\y-\B'\x)\x'
\end{equation}
vanishes at $\B=\B_0$. Then $\tita_0$ is a zero of the function $\Phi(\tita)=E_{H_0}\phi(\z;\tita)$.

By Lemma \ref{hipA6} there exists a bounded subset $\mathcal{C}$ and a function $\tita(\xi)$ such that $\tita_0$ is an interior point of $\tita(\mathcal{C})$ and since $\boldsymbol{\theta}_n\rightarrow\boldsymbol{\theta}_0$ a.s.,  $\tita_n\in \tita(\mathcal{C})$ for $n$ large enough, i.e., $\phi_{kj}(\z;\tita_n)$ and $\phi_{kj}(\z;\tita_0)$ belong to the Euclidean class $\mathfrak{F}_{kj}$ for $n$ sufficiently large. By (A5), the functions $\phi_{kj}(\z;\tita_n)$ and $\phi_{kj}(\z;\tita_0)$ are in the class $[\delta]$ of Lemma \ref{pollard} for each $\delta>0$ and $n$ sufficiently large. Hence, 
\begin{equation}
|\sqrt{n}\{\nu_n(\phi_{kj}(\cdot;\tita_n))-\nu_n(\phi_{kj}(\cdot;\tita_0))\}|\longrightarrow0
\end{equation}
in probability. 
Then since $\nu_n(\phi_{kj}(\cdot;\tita_n))-\nu_n(\phi_{kj}(\cdot;\tita_0))=\text{o}_P(1/\sqrt{n})$ for all $k=1,\hdots, q$ and $j=1,\hdots, p$ and $\phi_{kj}(\z;\tita)$ corresponds to the element $h=(j-1)q+k$ of the function $\phi$, we conclude that 
\begin{equation}\label{o chiquita}
\nu_n(\phi(\cdot;\tita_n))-\nu_n(\phi(\cdot;\tita_0))=\text{o}_P(1/\sqrt{n}).
\end{equation}

Since $\partial\Phi/\partial\ve{\B'}'$ is continuous in $\tita_0$, we have that 
\begin{equation}\label{expansion}
\Phi(\B,\Sig)=\Phi(\B_{0},\Sig)+
\left(\frac{\partial\Phi(\B,\Sig)}{\partial\ve{\B'}'}{(\B_{0},\Sig)}\right)\ve{\B'-\B_0'}+r(\tita)\ve{\B'-\B_0'}
\end{equation}
where $r(\tita)\rightarrow \mathbf{0}$ when $\tita\rightarrow \tita_0$.

Using a suitable change of variables, for all $\Sig\in \mathcal{S}_q$ we have that  
\begin{eqnarray*}
E_{H_0}W\left(d(\B_0,\Sig)\right)u_kx_j&=&(Ex_j)E_{F_0}W\left(\uu'\Sig^{-1}\uu\right)u_k\\
&=&(Ex_j)\left(\int_{\{\uu:u_k>0\}}u_kW\left(\uu'\Sig^{-1}\uu\right)f_0^*(\uu'\Sig_0^{-1}\uu)d\uu\right.\\
&+&\left.\int_{\{\uu:u_k<0\}}W\left(\uu'\Sig^{-1}\uu\right)u_kf_0^*(\uu'\Sig_0^{-1}\uu)d\uu\right)\\
&=&(Ex_j)\left(\int_{\{\uu:u_k>0\}}u_kW\left(\uu'\Sig^{-1}\uu\right)f_0^*(\uu'\Sig_0^{-1}\uu)d\uu\right.\\
&-& \left.\int_{\{\uu:u_k>0\}}u_kW\left(\uu'\Sig^{-1}\uu\right)f_0^*(\uu'\Sig_0^{-1}\uu)d\uu\right)=0.
\end{eqnarray*}
Since this holds for all $k=1,\hdots,q$ and $j=1,\hdots,p$ so that
\begin{equation}\label{Phi cero}
\Phi(\B_0,\Sig)=\mathbf{0}
\end{equation}
for all $\Sig\in \mathcal{S}_q$.

By \eqref{MMec1}, the pair $\tita_n=(\Bhat,\Sighat)$ is a zero of the function $(1/n)\sum_{i=1}^n\phi(\z_i,\tita)$. Using this, after doing some simple operations of sum and subtraction and using \eqref{o chiquita}, we have
\begin{eqnarray*}
	\mathbf{0}&=&(1/n)\sum_{i=1}^n\phi(\z_i,\tita_n)=	           E_{H_0}\phi(\z,\tita_n)+\left[\frac{1}{n}\sum_{i=1}^{n}\phi(\z_i,\tita_0)-E_{H_0}\phi(\z,\tita_0)\right]\\	&+&\left\{\frac{1}{n}\sum_{i=1}^{n}\left[\phi(\z_i,\tita_n)-\phi(\z_i,\tita_0)\right]-E_{H_0}\left[\phi(\z,\tita_n)-\phi(\z,\tita_0)\right]\right\}\\
&=&E_{H_0}\phi(\z,\tita_n)+\left[\frac{1}{n}\sum_{i=1}^{n}\phi(\z_i,\tita_0)-E_{H_0}\phi(\z,\tita_0)\right]\\
&+&\left[\nu_n(\phi(\cdot;\tita_n))-\nu_n(\phi(\cdot;\tita_0))\right]\\
&=& E_{H_0}\phi(\z,\tita_n)+\left[\frac{1}{n}\sum_{i=1}^{n}\phi(\z_i,\tita_0)-E_{H_0}\phi(\z,\tita_0)\right]+\text{o}_P(1/\sqrt{n}).
\end{eqnarray*}

Since $E_{H_0}\phi(\z,\tita_n)$ is equal to $\Phi(\Bhat,\Sighat)$, we can solve for $E_{H_0}\phi(\z,\tita_n)$ in the above equation and replace it in the expansion \eqref{expansion} for $\Phi(\Bhat,\Sighat)$, together with \eqref{Phi cero} we obtain the result
\begin{eqnarray*}
\mathbf{0}&=&\left(\frac{\partial\Phi(\B,\Sig)}{\partial\ve{\B'}'}{(\B_{0},\Sighat)}\right)\ve{\Bhat'-\B_0'}+r(\tita_n)\ve{\Bhat'-\B_0'}\\&+&\left[\frac{1}{n}\sum_{i=1}^{n}\phi(\z_i,\tita_0)-E_{H_0}\phi(\z,\tita_0)\right]+\text{o}_P({1}/{\sqrt{n}}).
\end{eqnarray*}
Since ${\partial\Phi}/{\partial\ve{\B'}'}$ is continuous in $\tita_0$ and as $r(\tita_n)=\text{o}_P(1)$, this reduces to
\begin{equation}\label{5.11}
\mathbf{0} = (\mathbf{\Lambda}+\text{o}_P(1))\ve{\Bhat'-\B_0'}+\left[\frac{1}{n}\sum_{i=1}^{n}\phi(\z_i,\tita_0)-E_{H_0}\phi(\z,\tita_0)\right]+\text{o}_P({1}/{\sqrt{n}}).
\end{equation}
According to the Central Limit Theorem \[\left[\frac{1}{n}\sum_{i=1}^{n}\phi(\z_i,\tita_0)-E_{H_0}\phi(\z,\tita_0)\right]=\text{O}_P({1}/{\sqrt{n}}),\] 
and since $\mathbf{\Lambda}$ is nonsingular, from \eqref{5.11} we get that $\ve{\Bhat'-\B_0'}=\text{O}_P({1}/{\sqrt{n}})$. 

Then \eqref{5.11} can be rewritten as 
\[
\mathbf{0}=\mathbf{\Lambda}\ve{\Bhat'-\B_0'}+\left[\frac{1}{n}\sum_{i=1}^{n}\phi(\z_i,\tita_0)-E_{H_0}\phi(\z,\tita_0)\right]+\text{o}_P({1}/{\sqrt{n}}).
\]

As we saw at the beginning of the proof, $\tita_0$ is a zero of $\Phi(\tita)=E_{H_0}\phi(\z;\tita)$, and therefore
\[
\sqrt{n}\ve{\Bhat'-\B_0'}=-\mathbf{\Lambda}^{-1}\frac{1}{\sqrt{n}}\sum_{i=1}^{n}\phi(\z_i,\tita_0)+\text{o}_P({1}).
\]

Since $\phi_{kj}(\z,\tita_0)$ has finite mean and covariance for each $k=1,\hdots, q$ and $j=1,\hdots, p$, the Theorem is proved after applying the Central Limit Theorem.\QED


\bigskip
\noindent\textbf{\underline{Proof of Proposition \ref{V}:}} Consider first the case $\Sig_0=\Id_q$.
The matrix $\mathbf{\Lambda}$ defined in \eqref{Lambda} can also be expressed as
\[
\mathbf{\Lambda}=\frac{\partial E_{H_0}W_1\left(d(\B,\Sig)^2\right)\ve{(\y-\B'\x)\x'}}{\partial \ve{\B'}'}{(\B_0,\sigma_0^2\Id_q)}.
\]
Since $W_1$ is differentiable with bounded derivative we can differentiating inside the expectation. 
We can now proceed analogously to the proof of \eqref{TFI2} and we have that
\begin{eqnarray*}
\mathbf{\Lambda}&=&
-\left[\frac{E_{F_0}2W_1'\left(\|\uu\|^2/{\sigma_0^2}\right)\|\uu\|^2}{q\sigma_0^2}+E_{F_0}W\left(\frac{\|\uu\|}{\sigma_0}\right)\right](E_{G_0}\x\x'\otimes\Id_q).
\end{eqnarray*}

Using the same arguments as before and $W(u)=\psi_1(u)/u$, we obtain
\begin{eqnarray*}
M&=&\left(E_{F_0}\left[W\left(\frac{\|\uu\|}{\sigma_0}\right)\right]^2\frac{\|\uu\|^2}{q}\right)E_{G_0}\x\x'\otimes\Id_q\\
&=&\left(E_{F_0}\left[\psi_1\left(\frac{\|\uu\|}{\sigma_0}\right)\right]^2\frac{\sigma_0^2}{q}\right)E_{G_0}\x\x'\otimes\Id_q.
\end{eqnarray*}

Since $(E_{G_0}\x\x'\otimes\Id_q)^{-1}=(E_{G_0}\x\x')^{-1}\otimes\Id_q$, the Proposition is proved for the case $\Sig_0=\Id_q$.

For the general case, let $\R$ a matrix such that $\Sig_0=\R\R'$ and consider the following transformation $\y^*=\R^{-1}\y$. Then $\B_0^*=\B_0\R'^{-1}$ and $\y^*={\B_0^{*}}'\x+\uu^*$, with $\uu^*=\R^{-1}\uu$. Observe that the distribution of $\uu^*$ is given by the density \eqref{densidad} with $\Sig_0=\Id_q$ and
\[
\ve{\B_0'}=\ve{\R{\B_0^*}'}=(\Id_p\otimes\R)\ve{{\B_0^*}'},
\]
and therefore, by the affine-equivariance of the MM-estimates, \eqref{V explicita} follows.\QED


\bigskip
\noindent\textbf{\underline{Proof of Theorem \ref{conv-algoritmo}:}}
We denote the weight $\omega_{ik}=W\left(d_i(\widetilde{\mathbf{B}}^{(k)},\widetilde{\mathbf{\Sig}}^{(k)})\right)$ by $\omega_i$ for each $1\leq i\leq n$ and $\tita^{(k)}=(\widetilde{\mathbf{B}}^{(k)},\widetilde{\mathbf{\Sig}}^{(k)})$ for each $k\geq 1$. Then, since $W(u)$ is nonincreasing in $|u|$ if and only if $\rho_1$ is concave (see page 326 of Maronna et al. \cite{Maronna}), we have
\begin{eqnarray}\label{ec1}
\nonumber\hspace{-1.9em}\sum_{i=1}^n\rho_1\left(d_i(\tita^{(k+1)})\right)-\hspace{-1.5em}&&\hspace{-0.5em}\sum_{i=1}^n\rho_1\left(d_i(\tita^{(k)})\right)\\
&&\hspace{-1.9em}\leq\frac{1}{2\sighat^2}\sum_{i=1}^n\omega_i \left[d_i^2(\widetilde{\mathbf{B}}^{(k+1)},\widetilde{\mathbf{\Ga}}^{(k+1)})-d_i^2(\widetilde{\mathbf{B}}^{(k)},\widetilde{\mathbf{\Ga}}^{(k)})\right]\hspace{-0.3em},
\end{eqnarray}
where $\widetilde{\mathbf{\Ga}}^{(k+1)}=\widetilde{\mathbf{\Sig}}^{(k+1)}/|\widetilde{\mathbf{\Sig}}^{(k+1)}|^{1/q}$ and $\widetilde{\mathbf{\Ga}}^{(k)}=\widetilde{\mathbf{\Sig}}^{(k)}/|\widetilde{\mathbf{\Sig}}^{(k)}|^{1/q}$.

Recall that for any positive definite matrix $\A$, the matrix $\widetilde{\mathbf{B}}^{(k+1)}$ minimizes 
\[
\sum_{i=1}^n\omega_id_i^2(\B,\A).
\]
Then
\[
\sum_{i=1}^n\omega_i d_i^2(\widetilde{\mathbf{B}}^{(k+1)},\widetilde{\mathbf{\Ga}}^{(k)})\leq \sum_{i=1}^n\omega_i d_i^2(\widetilde{\mathbf{B}}^{(k)},\widetilde{\mathbf{\Ga}}^{(k)})
\]
and therefore the sum on the right side of \eqref{ec1} is not greater than
\begin{equation}\label{ec2}
\sum_{i=1}^n\omega_i d_i^2(\widetilde{\mathbf{B}}^{(k+1)},\widetilde{\mathbf{\Ga}}^{(k+1)})- \sum_{i=1}^n\omega_i d_i^2(\widetilde{\mathbf{B}}^{(k+1)},\widetilde{\mathbf{\Ga}}^{(k)}).
\end{equation}
Since 
\[
\displaystyle{{\widetilde{\mathbf{\Ga}}^{(k+1)}}=\frac{\widetilde{\mathbf{C}}^{(k+1)}}{|\widetilde{\mathbf{C}}^{(k+1)}|^{1/q}}} \text{ con } \displaystyle{{\widetilde{\mathbf{C}}^{(k+1)}}=\sum_{i=1}^n\omega_{i}\uuhat_i(\widetilde{\mathbf{B}}^{(k+1)})\uuhat_i'(\widetilde{\mathbf{B}}^{(k+1)})}
\]
we have that $\widetilde{\mathbf{\Ga}}^{(k+1)}$ is the sample covariance matrix of the weighted residuals\\
\noindent$\sqrt{\omega_{i}} \uuhat_i(\widetilde{\mathbf{B}}^{(k+1)})$ normalized to unit determinant, which minimizes the sum of squared Mahalanobis norms of weighted residuals $\sqrt{\omega_{i}} \uuhat_i'(\widetilde{\mathbf{B}}^{(k+1)})$ among the matrices with determinant one, i.e., for any positive definite matrix $\V$ with $|\V|=1$
\[
\sum_{i=1}^n\omega_i d_i^2(\widetilde{\mathbf{B}}^{(k+1)},\widetilde{\mathbf{\Ga}}^{(k)})\leq \sum_{i=1}^n\omega_i d_i^2(\widetilde{\mathbf{B}}^{(k+1)},\V).
\]
Then, since $|\widetilde{\mathbf{\Ga}}^{(k+1)}|=|\widetilde{\mathbf{\Ga}}^{(k)}|=1$, we have that \eqref{ec2} is $\leq0$. \QED

\bigskip
\bigskip

\noindent{\bf Acknowledgements: \/}{We would like to thank the referees and the editor of the Journal of Multivariate Analysis for their helpful comments and suggestions. We also gratefully acknowledge the many important comments of Hendrik Lopuha\"a. This research was partially supported by grants PIP 216 from CONICET and PICT 899 from ANPCyT, Argentina.}

\normalsize{

}
\end{document}